\documentclass[11pt]{article}

\usepackage{amssymb}
\usepackage{amscd}
\usepackage{amsmath}
\usepackage{amsfonts}
\usepackage{theorem}
\usepackage{mathrsfs}
\usepackage{tikz}
\usepackage{graphicx}
\usepackage{float}
\usepackage{color}
\definecolor{red}{rgb}{.7,0,0} 
\definecolor{blue}{rgb}{0,0,1}

\def\bolita{\scriptscriptstyle\bullet}

{\theorembodyfont{\slshape}
\newtheorem{theorem}{Theorem}
\numberwithin{theorem}{section}
\newtheorem{proposition}[theorem]{Proposition}
\newtheorem{lemma}[theorem]{Lemma}
\newtheorem{corollary}[theorem]{Corollary}

}
{\theorembodyfont{\rmfamily}
\newtheorem{definition}[theorem]{Definition}

\newtheorem{remark}[theorem]{Remark}

}
\def\proof{{\noindent\sc Proof. \quad}}
\newcommand{\proofof}[1]{{\noindent\sc Proof of #1. \quad}}
\def\eproof{{\mbox{}\hfill\qed}\medskip}
\newcommand\qed{{\unskip\nobreak\hfil\penalty50\hskip2em\vadjust{}
\nobreak\hfil$\Box$\parfillskip=0pt\finalhyphendemerits=0\par}}

\def\R{{\mathbb{R}}}

\def\N{{\mathbb{N}}}
\def\C{{\mathbb{C}}}
\def\IS{{\mathbb{S}}}
\renewcommand{\P}{\mathbb{P}}
\def\E{\mathop{\mathbb{E}}}
\newcommand{\Exp}{\mathop{\mathbb{E}}}
\def\Prob{\mathop{\mathsf{Prob}}}

\def\mcU{\mathcal{U}}

\def\Cnn{{\C^{n\times n}}}

\def\Oh{{\mathcal{O}}}
\def\msD{\mathscr D}

\def\msU{\mathscr U}
\def\dist{\mathsf{dist}}
\def\diag{\mathsf{diag}}
\def\NJ{\mathrm{NJ}}

\def\cond{\mathsf{cond}}

\def\Id{\mathsf{Id}}

\def\uno{\mbox{1\hspace*{-2.5pt}l}}
\def\ndet{\mathrm{ndet}\,}

\renewcommand{\tilde}{\widetilde}
\def\a{\alpha}
\def\e{\varepsilon}

\def\d{\delta}
\def\s{\sigma}
\def\dS{d_{\IS}}
\def\dpr{d_{\P}}
\def\dR{d_{\mathsf{R}}}

\def\EC{{\sf Path-follow}}

\def\bH{K_{(\lambda,v)}}

\def\aV{\mathcal{\hat{V}}}

\def\oa{\overline{a}}
\def\oA{\overline{A}}
\def\oM{\overline{M}}
\def\oK{\overline{K}}

\def\oB{\overline{B}}
\def\oQ{\overline{Q}}
\def\Cnn{\C^{n\times n}}

\def\Dn{n}

\def\aS{\C\times \IS(\C^n)}

\def\av{\mathsf{av}}
\def\mum{\mu_{\max}}
\def\lin{\mathsf{lin}}

\newcommand{\algoritmo}{\begin{minipage}{0.87\hsize}\linea}
\newcommand{\falgoritmo}{\linea\end{minipage}\bigskip}
\newcommand{\linea}{\vspace*{-5pt}\hrule\vspace*{5pt}}
\newtheorem{algorithm}{Algorithm}
\def\espacio{\hspace*{1cm}}

\newcommand{\inputalg}[1]{\linea\bf Input:\quad\rm #1\vspace*{3pt}}
\newcommand{\specalg}[1]{\bf Preconditions:\quad\rm #1}
\newcommand{\Output}[1]{\linea\bf Output:\quad\rm #1\vspace*{2pt}}
\newcommand{\postcond}[1]{\bf Postconditions:\quad\rm #1\vspace*{3pt}}
\newcommand{\bodyalg}[1]{\linea\tt #1\vspace*{3pt}}

\def\la{\lambda}
\def\avcost{{\mathsf{Avg\_Cost}}}
\def\smcost{{\mathsf{Smoothed\_Cost}}}
\def\aviter{{\mathsf{Avg\_Num\_Iter}}}
\def\numiter{{\mathsf{Num\_Iter}}}
\def\siter{{\mathsf{Smoothed\_Num\_Iter}}}
\medskip
\def\bfuno{\mathbf{1}}

\def\mcN{\mathcal{N}}

\newcommand{\pes}[2]{\langle #1,#2\rangle}

\newcommand{\pcn}{\P(\C^n)}

\newcommand{\Un}{\mathcal{U}(n)}
\newcommand{\vol}{\mbox{vol}}

\newcommand{\ambient}{\Cnn\times\C\times\prc}
\newcommand{\ambientu}{\Cnn\times\C}
\newcommand{\prc}{\P(\mathbb{C}^n)}

\newcommand{\V}{\mathcal{V}}
\newcommand{\W}{\mathcal{W}}

\begin{document}

\bibliographystyle{plain}

\makeatletter


\def\JACM{Journal of the ACM}
\def\CACM{Communications of the ACM}
\def\ICALP{International Colloquium on Automata, Languages
            and Programming}
\def\STOC{annual ACM Symp. on the Theory
          of Computing}
\def\FOCS{annual IEEE Symp. on Foundations of Computer Science}
\def\SIAM{SIAM Journal on Computing}
\def\SIOPT{SIAM Journal on Optimization}
\def\MOR{Math. Oper. Res.}
\def\BSMF{Bulletin de la Soci\'et\'e Ma\-th\'e\-ma\-tique de France}
\def\CRAS{C. R. Acad. Sci. Paris}
\def\IPL{Information Processing Letters}
\def\TCS{Theoretical Computer Science}
\def\BAMS{Bulletin of the Amer. Math. Soc.}
\def\TAMS{Transactions of the Amer. Math. Soc.}
\def\PAMS{Proceedings of the Amer. Math. Soc.}
\def\JAMS{Journal of the Amer. Math. Soc.}
\def\LNM{Lect. Notes in Math.}
\def\LNCS{Lect. Notes in Comp. Sci.}
\def\JSL{Journal for Symbolic Logic}
\def\JSC{Journal of Symbolic Computation}
\def\JCSS{J. Comput. System Sci.}
\def\JoC{J. of Complexity}
\def\MP{Math. Program.}
\sloppy

\begin{title}
{{\bf  A stable, polynomial-time algorithm for the eigenpair problem}} 
\end{title}
\author{Peter B\"urgisser\thanks{Partially funded by DFG research grant BU 1371/2-2}\\
Technische Universit\"at Berlin\\
Institut f\"ur Mathematik\\
10623 Berlin, GERMANY\\
e-mail: {\tt pbuerg@math.tu-berlin.de}\\
\and
Felipe Cucker\thanks{Partially funded by
a GRF grant from the Research Grants Council of the
Hong Kong SAR (project number CityU 100810).}\\
Department of Mathematics\\
City University of Hong Kong\\
HONG KONG\\
e-mail: {\tt macucker@cityu.edu.hk}}

\date{}

\makeatletter
\maketitle
\makeatother

\thispagestyle{empty}

\begin{quote}
{\small 
{\bf Abstract.} 
We describe algorithms for computing eigenpairs 
(eigenvalue--eigenvector) of a complex $n\times n$ 
matrix $A$. These algorithms are numerically stable, 
strongly accurate, and theoretically efficient (i.e., polynomial-time). 
We do not believe they outperform in practice the algorithms 
currently used for this computational problem. The merit 
of our paper is to give a positive answer to a 
long-standing open problem in numerical linear algebra. 
}
\end{quote}
\medskip

\hfill\begin{minipage}{6cm}
{\small{\em  
So the problem of devising an algorithm [for the eigenvalue 
problem] that is numerically stable and globally (and quickly!) 
convergent remains open.}

\hfill J.~Demmel~\cite[page~139]{Demmel97}}
\end{minipage}
\bigskip\bigskip

\section{Introduction}

\subsection{The problem}

The quotation from Demmel opening this article, though 
possibly puzzling for those who day-to-day satisfactorily 
solve eigenvalue problems, summarizes a long-standing 
open problem in numerical linear algebra. The first algorithm 
that comes to mind for computing eigenvalues 
---to compute the characteristic polynomial 
$\chi_A$ of $A$ and then compute (i.e., approximate) its zeros---
has proved to be numerically unstable in practice. The so called 
Wilkinson's polynomial, 
$$
  w(x):=\prod_{i=1}^{20}(x-i) 
  = x^{20}+w_{19}x^{19}+\cdots+w_1x+w_0
$$
is often used to illustrate this fact. For a diagonal matrix $D$ 
with diagonal entries $1,2,\ldots,20$ (and therefore with 
$\chi_D(x)=w(x)$) an error of $2^{-23}$ in 
the computation of $w_{19}=210$ produces, even if the rest of the 
computation is done error-free, catastrophic variations in the zeros 
of $\chi_D$. For instance, the zeros at 18 and 19 collide into a double 
zero close to 18.62, which will unfold into two complex conjugate 
zeros if the error in $w_{19}$ is further increased. And yet, there 
is nothing wrong in the nature of $D$ (in numerical analysis terms, 
and we will be more detailed below, 
$D$ is a well-conditioned matrix for the eigenvalue problem). 
The trouble appears to lie in the method. 

Barred from using this immediate algorithm due to its numerical 
unstability, researchers devoted efforts to come up with alternate 
methods which would appear to be stable. Among those proposed, 
the one that is today's algorithm of choice is the iterated QR with 
shifts. This procedure behaves quite efficiently in general and yet, 
as Demmel pointed out in 1997~\cite[p.~139]{Demmel97}, 
\begin{quote}
{\small 
after more than 30 years of dependable service,
convergence failures of this algorithm have quite recently been
observed, analyzed, and patched [\dots]. But there is still no global
convergence proof, even though the current algorithm is considered
quite reliable. 
}
\end{quote}
Our initial quotation followed these words in Demmel's text. It 
demanded  
for an algorithm which will be numerically stable and for which, 
convergence, and if possible small complexity bounds, can be 
established. Today, 17 years after Demmel's text, this demand retains 
all of its urge. 

The only goal of this paper is to give a positive answer to it. 

\subsection{A few words on approximations}

It must be said upfront that we do not think the algorithm we 
propose will outperform, in general, iterated QR with shifts.  
It nonetheless possesses some worthy features which we 
want to describe in this introduction. The key one, we already 
mentioned, is that both convergence and complexity bounds 
can be established for it. It is also numerically stable. But in 
addition, it is {\em strongly accurate}. 

A starting point to understand the meaning of this last claim, 
is the observation that there are two different obstructions to 
the exact computation of an eigenvalue. Firstly, the use of finite 
precision, and the ensuing errors accumulating during the 
computational process. The expression {\em numerically stable} is
usually vested on algorithms for which this accumulated error on the 
computed quantities is not much 
larger than that resulting from an error-free computation on an 
input datum which has been read (and approximated) with machine 
precision. Secondly, the nonlinear character of the equations 
defining eigenvalues in terms of the given matrix. For $n\geq 5$, 
we learned from Galois, we cannot write down these eigenvalues 
in terms of the matrix' entries, not even using radicals. Hence, 
we can only compute approximations of them and this is so 
{\em even assuming infinite precision in the computation}. 

The expression {\em strongly accurate} refers to the quality of these 
approximations. It is common to find in the literature (at least) 
three notions of approximation which we next briefly 
describe. To simplify, we illustrate with the computation of 
a value $\zeta\in\C$ from a datum $A\in\C^N$ (and the reader 
may well suppose that this computation is done with infinite 
precision). We let $\tilde{\zeta}$ be the quantity actually computed
and we consider the following three requirements on it:
\begin{description}
\item{\em Backward approximation.}
The element $\tilde{\zeta}$ is the solution of a datum $\tilde{A}$ 
close to~$A$. Given $\e>0$, we say that $\tilde\zeta$ is an 
$\e$-backward approximation when $\|A-\tilde{A}\|\leq \e$. 
\item{\em Forward approximation.} 
The quantity $\tilde{\zeta}$ is close to $\zeta$. 
Given $\e>0$, we say that $\tilde\zeta$ is an 
$\e$-forward approximation when $|\zeta-\tilde{\zeta}|\leq \e$. 
\item{\em Approximation \`a la Smale.} An appropriate version of Newton's 
iteration, starting at $\tilde{\zeta}$, converges immediately, quadratically 
fast, to $\zeta$.  
\end{description}
These requirements are increasingly demanding. For instance, 
if $\zeta$ is an $\e$-backward approximation then the forward error 
$|\zeta-\tilde{\zeta}|$ is bounded, roughly, by $\e\,\cond(A)$.  
Here $\cond(A)$ is the condition number of $A$, 
a quantity usually greater than 1. So, in general, 
$\e$-backward approximations are not $\e$-forward approximations,  
and if $A$ is poorly conditioned $\tilde{\zeta}$ may be a 
poor forward approximation of $\zeta$. 
We also observe that if $\tilde{\zeta}$ is an approximation \`a la Smale 
we can obtain an $\e$-forward approximation from $\tilde{\zeta}$ 
by performing $\Oh(\log|\log\e|)$ Newton's steps. Obtaining an 
approximation \`a la Smale from an $\e$-forward approximation 
is a much less obvious process. 
\smallskip

When we say that our algorithm is strongly accurate, we refer to 
the fact that the returned eigenpairs are approximations \`a la Smale 
of true eigenpairs. 

\subsection{A few words on complexity}

The cost, understood as the number of arithmetic operations 
performed, of computing an approximation 
of an eigenpair for a matrix $A\in\Cnn$, depends on the 
matrix $A$ itself. Actually, and this is a common feature in 
numerical analysis, it depends on the condition $\cond(A)$ of 
the matrix $A$. But this condition is not known a priori. It was 
therefore advocated by Smale~\cite{Smale97} to eliminate this 
dependency in complexity bounds by endowing data space 
with a probability distribution and estimating 
expected costs. This idea has its roots in early work of 
Goldstine and von Neumann~\cite{vNGo51}. 

In our case, data space is $\Cnn$, and a common 
probability measure to endow it with is the standard 
Gaussian. Expectations of cost w.r.t.~this measure 
yield expressions in $n$ usually referred to as 
{\em average cost}. A number of considerations, including 
the suspicion that the use of the standard Gaussian could 
result in complexity bounds which are too optimistic compared 
with ``real life'', prompted Spielman and Teng to introduce 
a different form of probabilistic analysis, called {\em smoothed 
analysis}. In this, one replaces the average analysis' goal of 
showing that 
\begin{quote}
for a random $A$ it is unlikely that the cost for $A$ will be large
\end{quote}
by the following
\begin{quote}
for all $\oA$, it is unlikely that a slight random 
perturbation $A=\oA+\Delta A$ will require a large cost. 
\end{quote}
The expectations obtained for a smoothed analysis will now be 
functions of both the dimension $n$ and some measure of 
dispersion for the random perturbations (e.g., a variance). 

Smoothed analysis was first used for the simplex method of 
linear programming~\cite{ST:04}. Two survey expositions 
of its rationale are in~\cite{ST:02,ST:09}. One may argue that it 
has been well accepted by the scientific community from 
the fact that Spielman and Teng were awarded the 
G\"odel 2008 and Fulkerson 2009 prizes for it (the former 
by the theoretical computer science community and the latter 
by the optimization community). Also, in 2010, Spielman was awarded 
the Nevanlinna prize, and smoothed analysis appears in the 
laudatio of his work.  
\smallskip

In this paper we will exhibit bounds for the cost of our algorithm 
both for average and smoothed analyses. 

\subsection{A few words on numerical stability}

The algorithm we deal with in this paper belongs to the class 
of homotopy continuation methods. Experience has shown 
that algorithms in this class are very stable and stability 
analyses have been done for some of 
them, e.g.~\cite{BriCuPeRo}. Because of this, we will assume 
infinite precision all along this paper and steer clear of 
any form of stability analysis. We nonetheless observe that 
such an analysis can be easily carried out following the steps 
in~\cite{BriCuPeRo}. 

\subsection{Previous and related work}

Homotopy continuation methods go back, at least, to the 
work of Lahaye~\cite{Lahaye}. A detailed survey of their 
use to solve polynomial equations is in~\cite{Li:03}. 

In the early 1990s Shub and Smale set up a program to 
understand the cost of solving square systems of complex 
polynomial equations using homotopy methods. In a 
collection of articles~\cite{Bez1,Bez2,Bez3,Bez4,Bez5}, known 
as the {\em B\'ezout series}, they put in place many of the notions 
and techniques that occur in this article. The B\'ezout series 
did not, however, conclusively settle the understanding of the cost 
above, and in 1998 Smale proposed it as the 17th in his list 
of problems for the mathematicians of the 21st 
century~\cite{Smale98}. The problem is not yet fully solved but 
significant advances appear 
in~\cite{BePa:09,BePa:11,BuCu11}. 

The results in these papers cannot be directly used for 
the eigenpair problem since instances of the latter are ill-posed as  
polynomial systems. But the intervening ideas can be 
reshaped to attempt a tailor-made analysis for the eigenpair 
problem. A major step in this direction was done by 
Armentano in his PhD thesis (see~\cite{Armentano:13}), 
where a condition number $\mu$ for the eigenpair problem 
was defined and exhaustively studied. A further step was 
taken in~\cite{ArmCuc} where $\mu$ was used to analyze a 
randomized algorithm for the Hermitian eigenpair problem. 
A difference between our paper and 
both~\cite{Armentano:13} and~\cite{ArmCuc} is that in the latter 
the technical development binds inputs and outputs (eigenvalues) 
together. We have found more natural to uncouple them. 

Our paper follows this stream of research. 

\subsection{Structure of the exposition}

The remaining of this paper is divided into two parts. 
In the first one, Section~\ref{ss:GeoFrame} below, 
we introduce all the technical preliminaries, we describe 
with details the algorithms, and we state our main 
results (Theorems~\ref{thm:main} and~\ref{thm:main2}). 
The condition number~$\mu$, Newton's method,  
the notion of approximate eigenpair, and 
Gaussian distributions are among these 
technical preliminaries. The second part, which 
occupies us in Sections~\ref{sec:condition_property} 
to~\ref{sec:main_proof}, is devoted to proofs.

\section{Preliminaries, Basic Ideas, and Main Result}\label{ss:GeoFrame}

\subsection{Spaces and Metric Structures}\label{subsec:spaces}

Let $\Cnn$ be the set of $n\times n$ complex matrices. 
We endow this complex linear space with the restriction of the 
real part of the {\em Frobenius Hermitian product} $\pes~~_F$ given by  
$$
  \pes{A}{B}_F:=\mbox{trace }(B^*A)=\sum_{i,j=1}^n a_{ij}\,
  \overline {b_{ij}},
$$
where $A=(a_{ij})$ and $B=(a_{ij})$. The {\em Frobenius norm} $\|~\|_F$ 
on $\C^{n\times n}$ is the norm induced by $\pes~~_F$. 

On the product vector space $\ambientu$ we introduce the 
canonical Hermitian inner product structure and its 
associated norm structure and (Euclidean) distance. 

The space  $\C^n$ is equipped with the canonical Hermitian inner 
product $\pes~~$. We denote by $\P(\C^n)$ its associated 
projective space. This is a smooth manifold which carries a natural
Riemannian metric, namely, the real part of the \emph{Fubini-Study
metric} on $\pcn$. The Fubini-Study metric is the Hermitian structure
on $\P(\C^n)$ given in the following way: for $x\in\C^n$, 
\begin{equation}\label{eq:defTip}
   \pes{w}{w'}_x:=\frac{\pes{w}{w'}}{\|x\|^2},
\end{equation}
for all $w,\,w'$ in the Hermitian complement $x^\perp$ of $x$ in
$\C^n$. We denote by $d_\P$ the associated Riemannian 
distance.

The space $\ambient$ is endowed with the Riemannian product
structure. 
\smallskip

The Hermitian structure in the spaces $\Cnn$ and $\Cnn\times \C$ 
naturally endows them with a notion of {\em angle}. The Riemannian 
distances $d_{\mathbb{S}}$ and $\dpr$ 
in the unit sphere $\IS(\Cnn)$ and the projective space $\P(\C^n)$, 
respectively, are given precisely by the angle between its arguments. 

In addition, we will consider the following  function on 
$((\Cnn\setminus\{0\})\times\C \times \P(\C^n))^2$  
given by 
$$
   \dist((A,\la,v),(A',\la',v'))^2 := 
   \bigg\|\frac{A}{\|A\|_F}-\frac{A'}{\|A'\|_F}\bigg\|_F^2
   +\bigg|\frac{\la}{\|A\|_F}-\frac{\la'}{\|A'\|_F}\bigg|^2+\dpr(v,v')^2.
$$
When restricted to matrices with Frobenius norm~1 this defines a distance 
on $\IS(\Cnn)\times\C\times \P(\C^n)$. In general it is not, but it  
remains a convenient measure due to the 
scaling invariance of the eigenpair problem. The space 
$\IS(\Cnn)\times\C\times \P(\C^n)$ is naturally endowed with 
a Riemannian distance $\dR$ given by replacing 
$\|A-A'\|_F$ in the definition of $\dist$ by the angle $\dS(A,A')$. 
Since chords are smaller than their subtending angles we trivially 
have
\begin{equation}\label{eq:chord}
  \dist((A,\la,v),(A',\la',v')) \leq \dR((A,\la,v),(A',\la',v'))
\end{equation}
for all $(A,\la,v),(A',\la',v')\in\IS(\Cnn)\times\C\times \P(\C^n)$.

\subsection{The Varieties $\V$, $\W$, $\Sigma'$ and $\Sigma$}
\label{subsec:SolVar}

We define the \textit{solution variety} for the eigenpair 
problem as
$$ 
  \V := \left\{ (A, \lambda, v) \in   \Cnn \times \C
  \times \P ( \C^{n}): \ (A-\lambda \Id)v=0 \right\}.
$$

\begin{proposition}\label{prop:Vsmooth}
The solution variety $\V$ is a smooth submanifold of 
$\ambient$, of the same dimension 
as $\Cnn$. 
\end{proposition}

\proof
See~\cite[Proposition~2.2]{Armentano:13}.
\eproof

The set $\V$ inherits the Riemannian structure of the ambient space.
Associated to it there are natural projections: 
\smallskip

\addtocounter{equation}{1}
\begin{equation}\tag*{\raisebox{30pt}{(\theequation)}}
\begin{tikzpicture}
\path (-0.25,0.3) node[right]{$\V$};
\draw[->] (0.15,0) -- (1.15,-1) node[above =3mm]{$\pi_2$};
\draw[->] (-0.15,0) -- (-1.15,-1) node[above =3mm]{$\pi_1$};
\path (-1.9,-1.3) node[right]{$\Cnn$};
\path (0.5,-1.3) node[right]{$\C\times\P(\C^n)$.};
\end{tikzpicture}
\end{equation}
Because of Proposition~\ref{prop:Vsmooth}, 
the derivative $D\pi_1$ at $(A,\lambda,v)$ 
is a linear operator between spaces of equal dimension. 
The \emph{subvariety $\W$ of well-posed triples} is the subset of
triples $(A,\lambda,v)\in\V$ for which $D\pi_1 (A,\lambda,v)$ is an
isomorphism. In particular, when $(A,\lambda,v)\in\W$, the 
projection $\pi_1$ has a branch of its inverse (locally defined) 
taking $A\in\Cnn$ to $(A,\lambda,v)\in\V$. This branch of 
$\pi_1^{-1}$ is called the \emph{solution map} at $(\lambda,v)$. 

For $v\in\P(\C^n)$ we denote by $T_v$ 
the tangent space of $\P(\C^n)$ at $v$. We then have
$T_v :=\{x \in \C^n \mid \langle x,w\rangle =0 \}$ for 
any $w\in\C^n$ such that $v=[w]$. Let   
$P_{v^\perp}\colon\C^n\to T_v$ be the orthogonal projection. 
Given $(A,\lambda,v)\in\ambient$, we 
let $A_{\lambda,v}:T_v\to T_v$ be the linear operator 
given by 
\begin{equation}\label{eq:defAlv}
  A_{\lambda,v}:=P_{v^\perp} \circ (A-\lambda \Id)|_{T_v} .
\end{equation}
We will prove (in Proposition~\ref{pro:NJquot} below) 
that the set of well-posed triples is given by
\begin{equation}\label{eq:Alv}
    \W=\{(A,\lambda,v)\in\V:\;A_{\lambda,v}\,\mbox{ is invertible}\}
\end{equation}
(see also Lemma~2.7 in~\cite{Armentano:13}).

We define $\Sigma':=\V \setminus \W$ to be the variety of 
\emph{ill-posed triples}, and $\Sigma=\pi_1(\Sigma')\subset\Cnn$ the
\emph{discriminant variety}, i.e., the subset of {\em ill-posed inputs}.

\begin{remark}
From (\ref{eq:Alv}) it is clear that the subset $\Sigma\rq{}$ is the
set of triples $(A,\lambda,v)\in \V$ such that $\lambda$ is an
eigenvalue of $A$ of algebraic multiplicity at least~2. It follows
that $\Sigma$ is the set of matrices $A\in\Cnn$ with multiple
eigenvalues. In particular, when $A\in\Cnn\setminus\Sigma$, 
the eigenvalues of $A$ are pairwise different and 
$\pi_1^{-1}(A)$ is the set of triples
$(A,\lambda_1,v_1),\ldots,(A,\lambda_n,v_n)$, where 
$(\lambda_i,v_i)$, $i=1,\ldots,n$, are the eigenpairs of~$A$.
\end{remark}

\begin{proposition}\label{prop:codim}
The discriminant variety $\Sigma\subset\Cnn$ is a complex 
algebraic hypersurface. Consequently, for all $n\geq2$, 
we have $\dim_{\R}\Sigma=n^2-2$. 
\end{proposition}

\proof
See~\cite[Proposition~20.18]{Condition}.
\eproof

\medskip

\subsection{Unitary invariance}\label{sec:unitaryaction}

Let $\Un$ be the group of $n\times n$ unitary matrices.  The group
$\Un$ naturally acts on $\prc$ by $U([w]):=[Uw]$. 
In addition, $\Un$ acts on
$\Cnn$ by conjugation (i.e., $U\cdot A:=UAU^{-1}$), and on
$\ambientu$ by $U\cdot (A,\lambda):= (UAU^{-1},\lambda)$. 
These actions
define an action on the product space $\ambient$, namely,
\begin{equation}\label{eq:UactionPP}
    U\cdot(A,\lambda,v):= (UAU^{-1},\lambda,Uv).
\end{equation}

\begin{remark}
The varieties $\V$, $\W$, $\Sigma\rq{}$, and $\Sigma$,
are invariant under the action of $\Un$ (see~\cite{Armentano:13} for details).
\end{remark}

\subsection{Condition of a triple}

In a nutshell, condition numbers measure the worst possible 
output error resulting from a small perturbation on the input data. 
More formally, a condition number is the operator 
norm of a solution map such as the branches of $\pi^{-1}$ 
mentioned in~\S\ref{subsec:spaces} above, 
(see~\cite[\S14.1.2]{Condition} for a general exposition).  

In the case of the eigenpair problem, one can define two 
condition numbers for eigenvalue and eigenvector, respectively. 
Armentano has shown, however, that one can merge the two 
in a single one (see Section~3 in~\cite{Armentano:13} for details).    
Deviating slightly from~\cite{Armentano:13}, 
we define the \emph{condition number} of
$(A,\lambda,v)\in\W$ as
\begin{equation}\label{eq:defmu}
    \mu(A,\lambda,v) := \|A\|_F \|A_{\lambda,v}^{-1} \|,
\end{equation}
where $\|~\|$ is the operator norm. 

\begin{remark}\label{rem:muscaling}
The condition number $\mu$ is invariant under the action of the 
unitary group $\Un$, i.e.,
$\mu(UAU^{-1},\lambda,Uv)=\mu(A,\lambda,v)$ for all $U\in\Un$.
Also, it is easy to see, $\mu$ is scale 
invariant on the first two components. That is, 
$\mu(sA,s\lambda,v)=\mu(A,\lambda,v)$ for all nonzero real $s$.
\end{remark}

\begin{lemma}[Lemma~3.8 in~\cite{Armentano:13}]\label{le:lb_mu}
For $(A,\la,v)\in\V$ we 
have $\mu(A,\la,v) \ge \frac{1}{\sqrt{2}}$. \hfill $\Box$ 
\end{lemma}

The essence of condition numbers is that the measure how much may  
outputs vary when inputs are slightly perturbed. The following result, 
which we will prove in Section~\ref{sec:condition_property}, 
quantifies this property for $\mu$. 

\begin{proposition}\label{prop:dotzdotf-spherical} 
Let $\Gamma:[0,1]\to \V$, $\Gamma(t)=(A_t,\lambda_t,v_t)$ be a 
smooth curve such that $A_t$ 
lies in the unit sphere of $\Cnn$, for all $t$. 
Write $\mu_t:=\mu(\Gamma(t))$. Then we have, for all $t\in[0,1]$, 
$$
     |\dot{\lambda_t}|\leq \sqrt{1+\mu_t^2}\;\|\dot{A_t}\|,
    \qquad
    \|\dot{v_t}\|\leq \mu_t\;\|\dot{A_t}\|.
$$
In particular,
$$
  \big\|\dot{\Gamma}(t)\big\|\leq 
\sqrt{6}\;\mu_t\;\|\dot{A_t}\|.
$$
\end{proposition}

\begin{remark}
Since the property of $A_{\lambda,v}$ being invertible is Zariski open on
$\ambient$, the condition number $\mu$ can be naturally extended to 
a Zariski open neighborhood of $\W$ in $\ambient$. 
We will denote this extension 
also by $\mu$. In addition, when $A_{\lambda,v}$ is non-invertible 
we will let $\mu(A,\lambda,v):=\infty$, so that now $\mu$ is well-defined 
in all of $\ambient$. 
\end{remark}

Condition numbers themselves vary in a controlled manner. 
The following Lipschitz property makes this statement precise.

\begin{theorem}\label{th:lipsch}
Let $A,A'\in\C^{n\times n}$ be such that $\|A\|_F = \|A'\|_F = 1$, let 
$v,v'\in\C^n$ be nonzero, and let  $\la,\la'\in\C$ be such that 
$Av=\lambda v$.  
Suppose that
$$
 \mu(A,\la,v)\;\dist((A,\la,v),(A',\la',v'))\;\le\;  
 \frac{\e}{12.5}
$$
for $0<\e<0.37$. Then we have 
$$
 \frac{1}{1+\e}\,\mu (A,\lambda,v)\ \le\ \mu(A',\la',v') \ \le\  
 (1+\e)\, \mu(A,\la,v) .
$$
\end{theorem}

We give the proof of Theorem~\ref{th:lipsch} in 
Section~\ref{se:lipschitz}.
\medskip

Condition numbers are generally associated to input data. 
In the case of a problem with many possible solutions 
(of which returning an eigenpair of a given matrix is a 
clear case) one can derive the condition of a data from a notion 
of condition for each of these solutions. A discussion of 
this issue is given in~\cite[Section~6.8]{Condition}. For the 
purposes of this paper, we will be interested in two such derivations: 
the {\em maximum condition number} of $A$,
$$
    \mum(A):=\max_{j\leq n}\mu(A,\lambda_j,v_j),
$$
and the {\em mean square condition number} of $A$,
$$
   \mu_{\av}(A):=\left(\frac1n\sum_{j=1}^n
   \mu^2(A,\lambda_j,v_j)\right)^{\frac12}.
$$

We close this paragraph observing that restricted to  
the class of normal matrices, the condition number $\mu$  
admits the following elegant characterization.

\begin{lemma}[Lemma~3.12 in~\cite{Armentano:13}]\label{lem:char}
Let $A\in\Cnn\setminus\Sigma$ be normal, and let 
$(\lambda_1,v_1),\ldots,(\lambda_n,v_n)$ be its eigenpairs. Then
\begin{equation}\tag*{\qed}
\mu(A,\lambda_1,v_1) = 
\frac{\|A\|_F}{\min_{j=2,\ldots,n}|\lambda_j-\lambda_1|}
\end{equation}
\end{lemma}


\subsection{Newton's method and approximate eigenpairs}

For a nonzero matrix $A\in\Cnn$, we 
define the \emph{Newton map} associated to $A$, 
$$
    N_A:\C\times(\C^n\setminus\{0\})\to \C\times(\C^n\setminus\{0\}), 
$$
by $N_A(\lambda,v)=(\lambda-\dot\lambda,v-\dot v)$ where
\begin{equation*}
  \dot v  ={A_{\lambda,v}}^{-1}\,
  P_{v^\perp}(A-\lambda\,\Id)v, \qquad 
  \dot\lambda =\frac{\pes{(A-\lambda\,\Id)(v-\dot v)}{v}}{\pes vv}.
\end{equation*}
This map is defined for every 
$(\lambda,v)\in\C\times(\C^n\setminus\{0\})$ such that 
${A_{\lambda,v}:=P_{v^\perp}(A-\lambda\,\Id)|_{T_v}}$ is invertible. 
It was introduced in~\cite{Armentano:13} as the Newton operator 
associated to 
the evaluation map $(\lambda,v)\mapsto (A-\lambda\,\Id)v$ for a 
fixed~$A$. See Section~4 of~\cite{Armentano:13} for more details. 

\begin{definition}
Given $(A,\lambda,v)\in\W$ we say that 
$(\zeta,w)\in\C\times(\C^n\setminus\{0\})$ is an 
{\em approximate eigenpair} of $A$ with associated eigenpair $(\lambda,v)$ 
when for all $k\geq 1$ the $k$th iterate $N_A^k(\zeta,w)$ 
of the Newton map at $(\zeta,w)$ is well defined and satisfies 
$$
  \dist\big((A,N_A^k(\zeta,w)),(A,\lambda,v)\big)\leq
  \left(\frac12\right)^{2^k-1} \dist\big((A,\zeta,w),(A,\lambda,v)\big).
$$
\end{definition} 

 \begin{remark}
Note that, if $N_A(\zeta,w)=(\zeta',w')$ then 
$N_{sA}(s\zeta,\beta w)=(s\zeta', \beta w')$, for every 
$s\in\C\setminus\{0\}$ and $\beta\in\C\setminus\{0\}$. 
Hence, $N_A$ is correctly defined on $\C\times\P(\C^n)$ 
and the notion of approximate eigenpair scales correctly
in the sense that if $(\lambda,v)$ is an approximate eigenpair 
of $A$ then $(s\lambda,v)$ is an approximate eigenpair of 
$sA$, for all $s\in\C\setminus\{0\}$.
 \end{remark}

\begin{remark}
The notion of approximate solution as a point where Newton's method 
converges to a true solution immediately and quadratically fast was 
introduced by Steve Smale~\cite{Smale86}. It allows to elegantly talk
about polynomial time without dependencies on pre-established 
accuracies. In addition, these approximate solutions are ``good 
approximations'' (as mentioned in the statement of the Main Theorem) 
in a very strong sense. The distance to the exact solution dramatically 
decreases with a single iteration of Newton's method.     
\end{remark}

The following result estimates, in terms of the 
condition of an eigenpair, the radius of a ball of approximate 
eigenpairs associated to it. For a proof 
see~\cite[Theorem~5]{Armentano:13}.

\begin{theorem}\label{th15.1}
Let $A\in\Cnn$ with $\|A\|_F=1$ and  
$(\lambda,v),(\lambda_0,v_0)\in\C\times(\C^n\setminus\{0\})$. 
If $(\lambda,v)$ is a well-posed eigenpair of~$A$ and 
$$
    \dist\big((\lambda,v),(\lambda_0,v_0)\big)<
   \frac{c_0}{\mu(A,\lambda,v)}
$$ 
then $(\lambda_0,v_0)$ is an approximate eigenpair of $A$ with 
associated eigenpair $(\lambda,v)$. One may choose $c_0=0.2881$.
\eproof
\end{theorem}

\begin{remark}
We note that $N_A(\zeta,w)$ can be computed from the matrix 
$A$ and the pair $(\zeta,w)$ in $\Oh(n^3)$ operations, since 
the cost of this computation is dominated by that of 
inverting a matrix.
\end{remark}

\medskip

\subsection{Gaussian Measures on $\Cnn$}

Let $\sigma>0$. We say that the complex random variable 
$Z=X+\sqrt{-1}Y$ has distribution 
$\mcN_{\C}(0,\sigma^2)$ when the real part $X$ and the
imaginary part $Y$ are independent and identically distributed
(i.i.d.) drawn from $\mcN(0,\frac{\sigma^2}{2})$, i.e., they are 
Gaussian centered random variables with variance 
$\frac{\sigma^2}{2}$. 

If $Z\sim\mcN_{\C}(0,\sigma^2)$ then its
density $\varphi\colon\C\to\R$ with respect to the Lebesgue measure 
is given by 
\begin{equation*}
   \varphi(z):=\frac{1}{\pi\sigma^2}e^{-\frac{|z|^2}{\sigma^2}}. 
\end{equation*}

We will write $v\sim\mcN_{\C}(0,\bfuno_n\sigma^2)$ 
to indicate that the vector
$v\in\C^n$ is random with i.i.d.~coordinates drawn from
$\mcN_{\C}(0,\sigma^2)$. Also, we say that 
$A\in\Cnn$ is (isotropic) \emph{Gaussian}, 
if its entries are i.i.d. Gaussian random
variables. In this case we write 
$A\sim\mcN_{\C}(0,\sigma^2\Id_{n\times n})$, 
or simply $A\sim\mcN(0,\sigma^2\Id)$ 
(since we will only deal with square complex 
matrices here and in general the dimension $n$ 
will be clear from the context). 

If $\oA\in\Cnn$ and $G\sim\mcN(0,\sigma^2\Id)$,  
we say that the random matrix $A=G+\oA$ has the 
{\em Gaussian distribution centered at $\oA$}, and we write 
$A\sim\mcN(\oA,\sigma^2\Id)$. The density of this 
distribution is given by
$$
  \varphi_{n\times n}^{\oA,\sigma}(A):=\frac{1}{(\pi\sigma^2)^{n^2}}\,
  e^{-\frac{\|A-\oA\|_F^2}{\sigma^2}}. 
$$

Crucial in our development is the following result giving 
a bound on average condition for Gaussian matrices 
arbitrarily centered. We will prove 
it in Section~\ref{se:muave}. 

\begin{theorem}\label{th:mu2-bound}
For $\oQ\in\Cnn$ and $\s>0$ we have
$$
   \Exp_{Q\sim \mcN(\oQ,\s^2\Id)}
   \Big(\frac{\mu_{\av}^2(Q)}{\|Q\|^2} \Big)\ \le\ \frac{en}{2\s^2}.
$$
\end{theorem}

\begin{remark} 
\begin{description}
\item[(i)]
We note that no bound on the norm of $\oQ$ is required here.
Indeed, using $\mu_{\av}(s Q)=\mu_{\av}(Q)$,
it is easy to see that the assertion for  a pair $(\oQ,\s)$ implies the
assertion for $(s\oQ,s\s)$, for any $s>0$.
\item[(ii)]
Because of Proposition~\ref{prop:codim}, with probability one, 
matrices drawn from $\mcN(\oQ,\s^2\Id)$ have all its
eigenvalues different. Therefore the expected value in 
Theorem~\ref{th:mu2-bound} is well-defined.  
\end{description}
\end{remark}

\subsection{Truncated Gaussians and smoothed analysis}

For $T,\sigma>0$, we define the 
{\em truncated Gaussian} $\mcN_T(0,\sigma^2\Id)$
on $\Cnn$ to be the distribution  
given by the density 
\begin{equation}\label{eq:truncated}
   \rho^{\sigma}_T(A)=
    \left\{\begin{array}{ll}
     \frac{\varphi_{n\times n}^{0,\sigma}(A)}{P_{T,\s}} & 
     \mbox{if $\|A\|_F\leq T$,}\\
     0 &\mbox{otherwise,}\end{array}\right.
\end{equation}
where $P_{T,\s}:=\Prob_{f\sim \mcN(0,\sigma^2\Id)}\{\|f\|\leq T\}$, 
and, we recall, $\varphi_{n\times n}^{0,\sigma}$ is the density of 
$\mcN(0,\sigma^2\Id)$. For the rest of this paper we fix 
the threshold $T:=\sqrt{2}\,n$. The fact that 
$\|A\|_F^2$ is chi-square distributed with $2n^2$ degrees of freedom, 
along with~\cite[Corollary~6]{choi:94} yield the following result.

\begin{lemma}\label{lem:X}
We have $P_{T,\s} \ge \frac12$
for all $0<\sigma\leq 1$.\eproof
\end{lemma}

The space $\C^{n\times n}$ of matrices with the Frobenius norm 
and the space $\C^{n^2}$ with the canonical Hermitian product are 
isomorphic as Hermitian product spaces. Hence, the 
Gaussian $\mcN(0,\sigma^2\Id_{n\times n})$ 
on the former corresponds to the 
Gaussian $\mcN(0,\sigma^2\bfuno_{n^2})$ 
on the latter, and we can deduce invariance 
of $\mcN(0,\sigma^2\Id)$ under the action of $\mcU(n^2)$ 
(in addition to that for conjugation under $\mcU(n)$ discussed 
in~\S\ref{sec:unitaryaction}). The same is true for the truncated 
Gaussian $\mcN_T(0,\sigma^2\Id)$. In particular, the pushforward 
of both distributions for the projection 
$\Cnn\setminus\{0\}\to\IS(\Cnn)$, $A\mapsto\frac{A}{\|A\|_F}$, is 
the uniform distribution $\msU(\IS(\Cnn))$ 
(see~\cite[Chapter~2]{Condition} for details) and we have
\begin{equation}\label{eq:truncating}
  \Exp_{A\sim\mcN(0,\sigma^2\Id)} F(A) \,=\,
  \Exp_{A\sim\mcN_T(0,\sigma^2\Id)} F(A) \,=\,
  \Exp_{A\sim\msU(\IS(\Cnn))} F(A).
\end{equation}
for any integrable scale invariant function $F:\Cnn\to\R$.
\medskip

Complexity analysis has traditionally been carried out either 
in the {\em worst-case} or in an {\em average-case}. More generally, 
for a function $F:\R^m\to\R_+$, the former amounts to the evaluation 
of $\sup_{a\in\R^m}F(a)$ and the latter to that of 
$\Exp_{a\sim\msD}F(a)$ for some probability distribution 
$\msD$ on $\R^m$. Usually, $\msD$ is taken to be 
the standard Gaussian $\mcN(0,\Id)$. With the beginning of 
the century, Daniel Spielman and Shang-Hua Teng introduced 
a third form of analysis, {\em smoothed analysis}, which is meant to 
interpolate between worst-case and average-case. We won't 
elaborate here on the virtues of smoothed analysis; a defense 
of these virtues can be found, e.g., 
in~\cite{ST:02,ST:09} or in~\cite[\S2.2.7]{Condition}. We 
will instead limit ourselves to the description of what smoothed 
analysis is and which form it will take in this paper. 

The idea is to replace the two operators above (supremum 
and expectation) by a combination of the two, namely,
$$
   \sup_{\oa\in\R^m} \Exp_{a\sim\msD(\oa,\sigma)} F(a) 
$$
where $\msD(\oa,\sigma)$ is a distribution ``centered'' at 
$\oa$ having $\sigma$ as a measure of dispersion. A  
typical example is the Gaussian $\mcN(\oa,\sigma^2\Id)$. 
Another example, used for scale invariant functions $F$, 
is the uniform measure on a spherical cap centered at 
$\oa$ and with angular radius $\sigma$ on the unit 
sphere $\IS(\R^m)$ (reference~\cite{Condition} exhibits  
smoothed analyses for both examples of distribution). 
In this paper we will perform a smoothed analysis with respect 
to a truncated Gaussian. More precisely, we will be interested in 
quantities of the form 
\begin{equation*}
   \sup_{\oA\in\Cnn} \Exp_{A\sim\mcN_T(\oA,\sigma^2\Id)} F(A) 
\end{equation*}
where $F$ will be a measure of computational cost for the 
eigenpair problem. 
We note that, in addition to the usual dependence on $n$, 
this quantity depends also on~$\sigma$ and tends to $\infty$ 
when $\sigma$ tends to $0$. When $F$ is scale invariant, 
as in the case of $\mu_{\av}$ or $\mum$, it is customary 
to restrict attention to matrices of norm 1. That is, to 
study the following quantity:
\begin{equation}\label{eq:smoothed}
   \sup_{\oA\in\IS(\Cnn)} \Exp_{A\sim\mcN_T(\oA,\sigma^2\Id)} F(A).  
\end{equation}

\subsection{The eigenpair  continuation algorithm}
\label{subsec:ALH}

We are ready to describe the main algorithmic construct 
in this paper. For the algorithmic purposes, it will be more convenient 
to view the solution variety as the corresponding subset of
$\Cnn\times\C\times(\C^n\setminus\{0\})$, which, abusing notation, 
we still denote by $\V$. 

Suppose that we are given an input matrix $A\in\Cnn$ and 
an \emph{initial triple} $(M,\lambda,v)$ in the solution 
variety~$\V\subseteq\Cnn\times\C\times(\C^n\setminus\{0\})$ 
such that $A$ and $M$ are $\R$-linearly independent. 
Let $\a:=\dS(M,A)\in(0,\pi)$ denote the {\em angle} between the 
rays $\R_+{A}$ and $\R_+{M}$. 
Consider the line segment
$[M,A]$ in $\Cnn$ with endpoints~$M$ and~$A$. 

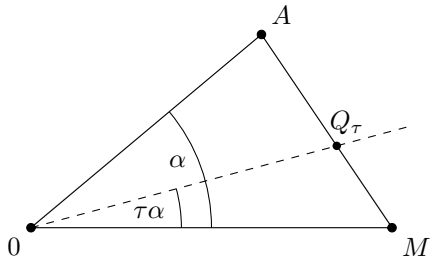
\begin{figure}[H]
\begin{center}
\begin{tikzpicture}[scale=4]
\def\mypoints{0.4}
\def\myalp{40};
\def\mybet{15};
\pgfmathsetmacro\sinalp{sin(\myalp)}
\pgfmathsetmacro\cosalp{cos(\myalp)}
\pgfmathsetmacro\sinbet{sin(\mybet)}
\pgfmathsetmacro\cosbet{cos(\mybet)}
\filldraw (0,0) circle(\mypoints pt) node[below left]{\small $0$} -- (1.2,0) circle(\mypoints pt) node[below right]{\small $M$} -- (\cosalp,\sinalp) circle(\mypoints pt) node[above right]{\small $A$} -- (0,0);
\draw (0.6,0) arc (0:\myalp:0.6);
\draw (0.5,0) arc (0:\mybet:0.5);
\draw[dashed] (0,0) -- +(\mybet:1.3cm);
\pgfpathcircle{
  \pgfpointintersectionoflines{\pgfpointxy{0}{0}}{\pgfpointxy{\cosbet}{\sinbet}}  {\pgfpointxy{1.2}{0}}{\pgfpointxy{\cosalp}{\sinalp}}
}{\mypoints pt}
\pgfusepath{fill}
\path (1.05,0.35) node{\small $Q_\tau$}
      (0.55,0.22) node[left]{\small $\alpha$}
      (0.48,0.05) node[left]{\small $\tau\alpha$};
\end{tikzpicture}
\end{center}
  \caption{The family $Q_\tau$, $\tau\in[0,1]$.}
  \label{fig:situation}
\end{figure}

We parameterize this segment by writing
$$
  [M,A]=\{Q_\tau\in\Cnn\mid \tau\in[0,1]\}
$$
with $Q_\tau$ being the only point in $[M,A]$ such that
$\dS(M,Q_\tau)=\tau\a$ (see Figure~\ref{fig:situation}).   

If the line segment $[M,A]$ does not intersect the 
discriminant variety~$\Sigma$, 
then starting at the eigenpair~$(\lambda,v)$ of $M$,  
the map $[0,1]\to \Cnn,\tau \mapsto Q_\tau$, 
can be uniquely extended to a continuous map
\begin{equation}\label{eq:curva}
  [0,1]\to \V,\quad \tau \mapsto (Q_\tau,\lambda_\tau,v_\tau), 
\end{equation}
such that $(\lambda_0,v_0)= (\lambda,v)$. 
We call this map the {\em lifting} of $[M,A]$ with origin $(M,\lambda,v)$.
We shall also call $\tau\mapsto (Q_\tau,\lambda_\tau,v_\tau)$ 
the {\em solution path} in~$\V$ corresponding to the input matrix~$A$ 
and initial triple $(M,\lambda,v)$. 

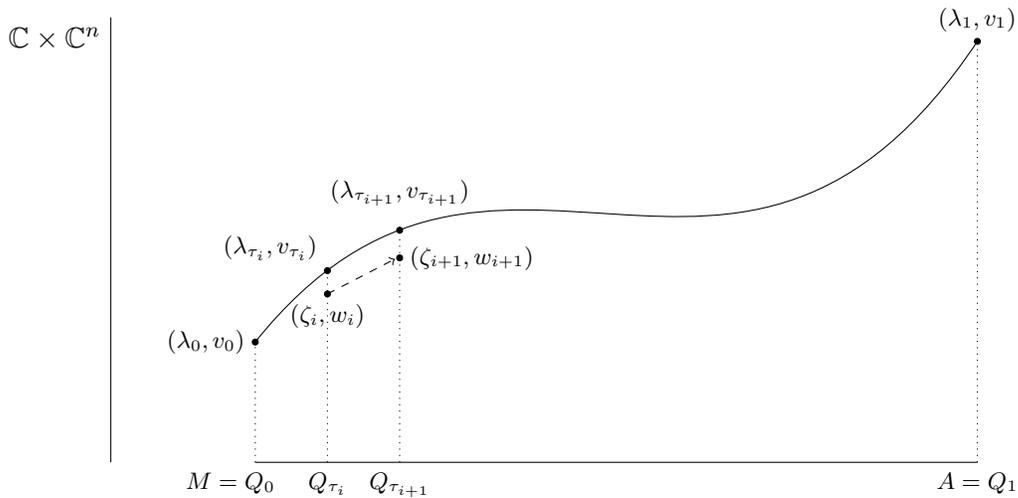
\begin{figure}[H]
\begin{center}
\begin{tikzpicture}[scale=1.6]
\def\mypoints{0.7}
\filldraw (0,0) circle(\mypoints pt); 
\filldraw (6,2.5) circle(\mypoints pt); 
\draw (0,0) node[left]{\footnotesize $(\lambda_0,v_0)$} .. controls (2,2.5) and (4,-0.5) .. (6,2.5) node[above]{\footnotesize $(\lambda_1,v_1)$};
\draw (0,-1)  node[below=2.8mm, left=-4mm]{\footnotesize $M=Q_0$} -- (6,-1)  node[below]{\footnotesize $A=Q_1$};
\draw[dotted] (0,-1) -- (0,0);
\draw[dotted] (0.6,-1)  node[below]{\footnotesize $Q_{\tau_i}$}-- (0.6,0.6) node[above left]{\footnotesize $(\lambda_{\tau_i},v_{\tau_i})$};
\filldraw (0.6,0.4) circle(\mypoints pt)  node[below]{\footnotesize $(\zeta_i,w_i)$}; 
\filldraw (0.6,0.594) circle(\mypoints pt); 
\draw[dotted] (1.2,-1)  node[below]{\footnotesize $Q_{\tau_{i+1}}$} -- (1.2,0.9); 
\draw (1.2,1.05) node[above]{\footnotesize $(\lambda_{\tau_{i+1}},v_{\tau_{i+1}})$};
\filldraw (1.2,0.7) circle(\mypoints pt) node[right]{\footnotesize $(\zeta_{i+1},w_{i+1})$}; 
\filldraw (1.2,0.93) circle(\mypoints pt); 
\draw[dotted] (6,-1) -- (6,2.5);
\draw[dashed] [->](0.6,0.4) -- (1.16,0.69);
\draw (-1.2, -1) -- (-1.2,2.7) node[below left]{$\C\times\C^n$};
\end{tikzpicture}
\end{center}
  \caption{The continuation of the solution path.}
  \label{fig:homotopy}
\end{figure}

Our algorithm relies on the obvious fact that the pair 
$(\lambda_1,v_1)$, corresponding to 
$\tau=1$, is an eigenpair of $A$. We therefore want to find 
an approximation of this pair and to do so, the idea is 
to start with the eigenpair $(\lambda,v)=(\lambda_0,v_0)$
of $M=Q_0$ and numerically follow the path 
$(Q_\tau,\lambda_\tau,v_\tau)$. This is done by 
subdividing the interval $[0,1]$ into subintervals with extremities 
at $0=\tau_0<\tau_1<\cdots<\tau_K=1$
and by successively computing approximations~$(\zeta_i,w_i)$ 
of~$(\lambda_{\tau_i},v_{\tau_i})$ by Newton's method. 
To ensure that these are good approximations, we actually  
want to ensure that for all $i$, $(\zeta_i,w_i)$ is 
an approximate eigenpair of $Q_{\tau_{i+1}}$.
Figure~\ref{fig:homotopy} attempts to convey the general idea. 

The following pseudocode gives a precise description of how 
this is done. The letter $\xi$ denotes a constant, namely 
$\xi=0.001461$. 
\bigskip\bigskip

\algoritmo
\begin{algorithm}\label{alg:EC}
{\sf Path-follow}\\
\inputalg{$A,M\in \Cnn$ and $(\lambda,v)\in\C\times\C^n$}\\
\specalg{$(M-\lambda\,\Id)v=0$, $M\not\in{\R}A$, $v\neq 0$}\\
\bodyalg{
$\a:=\dS(M,A)$, $r:=\|A\|_F$, $s:=\|M\|_F$\\[2pt]
$\tau:=0$, $Q:=M$, $(\zeta,w):=(\lambda,v)$\\[2pt]
repeat\\[2pt]
\espacio $\Delta\tau:= \frac{\xi}{\alpha 
\mu^2(Q,\zeta,w)}$\\[2pt]
\espacio $\tau:=\min\{1,\tau+\Delta\tau\}$\\[2pt]
\espacio $t:=\frac{s}{s+r(\sin\a\cot(\tau\alpha) - \cos\alpha)}$\\[2pt]
\espacio $Q:=tA+(1-t)M$\\[2pt]
\espacio $(\zeta,w):=N_{Q}(\zeta,w)$\\[2pt]
until $\tau= 1$\\[2pt]
return $(\zeta,w)$\\
}
\Output{$(\zeta,w)\in \C\times \C^{n}$}\\
\postcond{The algorithm halts if the lifting of 
$[M,A]$ at $(\lambda,v)$ does not cut $\Sigma'$. 
In this case, $(\zeta,w)$ is an approximate 
eigenpair of $A$.}
\end{algorithm}
\falgoritmo

Our next result estimates the number of iterations 
performed by algorithm~\EC. We summarize its proof in 
Section~\ref{sec:homotopy}.

\begin{proposition}\label{thm:main_path_following}
Suppose that $[M,A]$ does not intersect the discriminant
variety~$\Sigma$. 
Then the algorithm \EC\ stops after at most $K:=K(A,M,\lambda)$ 
steps with
$$
     K \ \le\  1077\,\,\dS(M,A)\int_0^1
     \mu^2(Q_\tau,\lambda_\tau,v_\tau)\,d\tau.
$$
The returned pair $(\zeta,w)$ is an approximate eigenpair of $A$ with
associated eigenpair $(\lambda_1,v_1)$. 
Furthermore, the bound above is optimal up to a constant: 
we have
$$
     K \ \ge\ 434\,\,\dS(M,A)\int_0^1
     \mu^2(Q_\tau,\lambda_\tau,v_\tau)\,d\tau.
$$
\end{proposition}

\subsection{Initial triples and global algorithms}

The \EC\ routine assumes an initial triple $(M,\lambda,v)$ 
at hand. We next describe a construction for such initial triples. 

For $k\in\N$ we consider the set of points  
$$
   S_k=\left\{\left(-1+\frac{2p}{k},-1+\frac{2q}{k}\right)\mid 
    0\leq p,q\leq k\right\}\subset\C.
$$
This is a set of $(k+1)^2$ points which are equidistributed 
on the square of side length~2, inscribed in the circle 
$\{z\in\C\mid |z|\leq \sqrt{2}\}$. 

\begin{figure}[H]
\begin{center}
\begin{tikzpicture}[scale=0.8]
\draw (0,0) circle (2.82842cm);
\draw[->] (0,-3.2) -- (0,3.2); 
\draw[->] (-3.2,0) -- (3.2,0); 
\path (-2,-2) node {$\bolita$};
\path (-1,-2) node {$\bolita$};
\path (-0,-2) node {$\bolita$};
\path (1,-2) node {$\bolita$};
\path (2,-2) node {$\bolita$};
\path (-2,-1) node {$\bolita$};
\path (-1,-1) node {$\bolita$};
\path (-0,-1) node {$\bolita$};
\path (1,-1) node {$\bolita$};
\path (2,-1) node {$\bolita$};
\path (-2,0) node {$\bolita$};
\path (-1,0) node {$\bolita$};
\path (-0,0) node {$\bolita$};
\path (1,0) node {$\bolita$};
\path (2,0) node {$\bolita$};
\path (-2,1) node {$\bolita$};
\path (-1,1) node {$\bolita$};
\path (-0,1) node {$\bolita$};
\path (1,1) node {$\bolita$};
\path (2,1) node {$\bolita$};
\path (-2,2) node {$\bolita$};
\path (-1,2) node {$\bolita$};
\path (-0,2) node {$\bolita$};
\path (1,2) node {$\bolita$};
\path (2,2) node {$\bolita$};
\end{tikzpicture}
\end{center}
  \caption{The set $S_2$.}
  \label{fig:initial}
\end{figure}
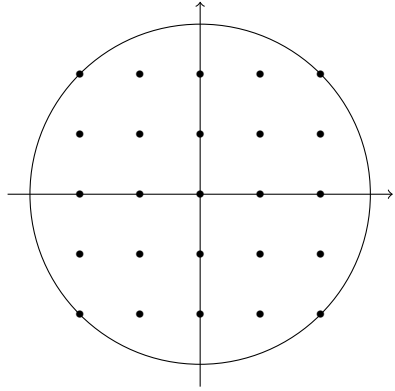

For $n\geq 2$ we let $k:=\lceil \sqrt{n}\rceil-1$ and 
define the matrix $D_n$ to be the diagonal matrix 
whose diagonal entries $z_1,\ldots,z_n$ are the 
first $n$ elements in $S_k$ (say for the lexicographical ordering). 
The eigenpairs of $D_n$ are the pairs $(z_j,e_j)$ for 
$j=1,\ldots,n$. 

\begin{lemma}\label{lem:initial}
For all $j\leq n$ we have $\mu(D_n,z_j,e_j)=\Theta(n)$.
\end{lemma}

\proof
Let $k:=\lceil \sqrt{n}\rceil-1$. Assume first that $k$ is even 
and write $z_{pq}:=\Big(-1+\frac{2p}{k},-1+\frac{2q}{k}\Big)$. Then
\begin{eqnarray*}
\sum_{p,q=0}^{k} |z_{pq}|^2 &\leq& 8 \sum_{p=0}^{\frac{k}{2}} 
 \sum_{q=0}^p  |z_{pq}|^2\; =\; 
8 \sum_{p=0}^{\frac{k}{2}} \sum_{q=0}^p \frac{p^2+q^2}{k^2}\\
&=& \frac{8}{k^2} \sum_{p=0}^{\frac{k}{2}}  
\left(p^2(p+1)+\sum_{q=0}^p q^2\right) 
\;=\; \frac{8}{k^2} \sum_{p=0}^{\frac{k}{2}}  \left(p^2(p+1)+ 
\frac{p^3}{3}+\frac{p^2}{2}+\frac{p}6\right)\\
&=& \frac{8}{k^2} \sum_{p=0}^{\frac{k}{2}}  
\left(\frac{4p^3}{3}+\frac{3p^2}{2}+\frac{p}6\right)
\;=\; \Theta(k^2) \;=\; \Theta(n). 
\end{eqnarray*}
One can similarly prove that the same bound holds for $k$ odd. 

It follows that for $n\geq 2$
$$
  \|D_n\|_F =\bigg(\sum_{p,q=0}^{\lceil \sqrt{n}\rceil-1} 
    |z_{pq}|^2\bigg)^{1/2} =\Theta(\sqrt{n}).
$$
But clearly, the smallest distance between two different eigenvalues 
of $D_n$ is $\frac{2}{k}=\Theta({n^{-\frac12}})$. 
And $D_n$ is diagonal, hence normal.
We conclude with Lemma~\ref{lem:char} 
that $\mu(D_n) =\Theta(n)$
as claimed. 
\eproof

We now put together the continuation algorithm and a 
specific initial triple. 
\bigskip\bigskip

\algoritmo
\begin{algorithm}\label{alg:LV}
{\sf Single\_Eigenpair}\\
\inputalg{$A\in \Cnn$}\\
\bodyalg{
compute $D_n$ \\[2pt]
set $M:=\frac{D_n}{\|D_n\|_F}$ \\[2pt]
randomly choose $j\in\{1,\ldots,n\}$ \\[2pt]
$(\zeta,w):=\mbox{\EC}(A,M,m_{jj},e_j)$\\[2pt]
}
\Output{$(\zeta,w)\in \C\times \C^{n}$}\\
\postcond{The algorithm halts if $[M,A]\cap\Sigma=\emptyset$. 
In this case, the pair $(\zeta,w)$ is an approximate 
eigenpair of $A$.}
\end{algorithm}
\falgoritmo

\begin{remark}\label{rem:twoThings}
The fact that the real codimension of $\Sigma$ in $\Cnn$ is 2 
(shown in Proposition~\ref{prop:codim}) ensures that, 
almost surely, the segment $[M,A]$ does not intersect 
$\Sigma$ and therefore, that almost surely {\sf Single\_Eigenpair} 
halts.  
\end{remark}

Given a matrix $A\in\Cnn$, the cost of 
algorithm~{\sf Single\_Eigenpair} 
with input $A$ depends on the triple 
$(M,m_{jj},e_j)$ which is random in the measure that $j$ is. 
We therefore consider the \emph{randomized cost} of this 
algorithm on input $A$. This amounts to 
the expected number of iterations of 
algorithm~{\sf Path-follow} with input $A$ times 
the $\Oh(n^3)$ cost of each iteration. The former 
is given by
$$
     \numiter(A):=\frac1n\sum_{j=1}^n K(A,M,m_{jj},e_j).
$$
Since we are interested in the \emph{average complexity} of 
{\sf Single\_Eigenpair} we will further take the expectation of 
$\numiter(A)$ when $A$ is drawn from $\mcN(0,\sigma^2\Id)$. 
We therefore consider  
$$
  \aviter(n):=\E_{A\sim\mcN(0,\sigma^2\Id)}\frac1n
   \sum_{j=1}^n K(A,M,m_{jj},e_j).
$$
Multiplying this expression by the cost $\Oh(n^3)$ of each iteration, 
we obtain the \emph{average cost} $\avcost(n)$ of {\sf Single\_Eigenpair}. 
\smallskip

We can also consider the smoothed cost of {\sf Single\_Eigenpair} 
by drawing the input matrix $A$ from $\mcN_T(\oA,\sigma^2\Id)$ 
where $\oA\in\IS(\Cnn)$ is arbitrary. We thus define
$$
  \siter(n):=\sup_{\oA\in\IS(\Cnn)}
   \E_{A\sim\mcN_T(\oA,\sigma^2\Id)}\frac1n
   \sum_{j=1}^n K(A,M,m_{jj},e_j)
$$
and multiplying by $\Oh(n^3)$ we obtain a corresponding notion 
of smoothed cost $\smcost(n)$. 

We can now state our first main result. 

\begin{theorem}\label{thm:main}
Algorithm~{\sf Single\_Eigenpair} returns (almost surely) an 
approximate eigenpair of its input $A\in\Cnn$. Its average 
cost satisfies
$$
    \avcost(n)=\Oh(n^8). 
$$ 
For every $0<\sigma\leq 1$, its smoothed cost satisfies
$$
    \smcost(n)=\Oh\Big(\frac{n^8}{\sigma^2}\Big). 
$$ 
\end{theorem} 

We can easily modify algorithm~{\sf Single\_Eigenpair} 
to compute all the eigenpairs.
\bigskip\bigskip

\algoritmo
\begin{algorithm}\label{alg:All}
{\sf All\_Eigenpairs}\\
\inputalg{$A\in \Cnn$}\\
\bodyalg{
compute $D_n$ \\[2pt]
set $M:=\frac{D_n}{\|D_n\|_F}$ \\[2pt]
for $j\in\{1,\ldots,n\}$ do\\[2pt]
\espacio $(\zeta_j,w_j):=\mbox{\EC}(A,M,m_{jj},e_j)$\\[2pt]
}
\Output{$\{(\zeta_1,w_1),\ldots,(\zeta_n,w_n)\}\in (\C\times \C^{n})^n$}\\
\postcond{The algorithm halts if $[M,A]\cap\Sigma=\emptyset$. 
In this case, the pairs $(\zeta_j,w_j)$ are approximate 
eigenpairs of $A$ with pairwise different associated eigenpairs.}
\end{algorithm}
\falgoritmo

This is no longer a randomized algorithm. In particular, the 
number of iterations performed by {\sf All\_Eigenpairs} 
on input $A$, which is now
$$
     \numiter(A):=\sum_{j=1}^n K(A,M,m_{jj},e_j),
$$
is no longer a random variable. We derive from these 
quantities the corresponding notions of 
$\avcost(n)$ and $\smcost(n)$, for which we state our second 
main result.

\begin{theorem}\label{thm:main2}
Algorithm~{\sf All\_Eigenpairs} returns (almost surely)
$n$~approximate eigenpairs of its input $A\in\Cnn$, 
with pairwise different associate eigenpairs. Its average 
cost satisfies
$$
    \avcost(n)=\Oh(n^9). 
$$ 
For every $\sigma\leq 1$ its smoothed cost satisfies
$$
    \smcost(n)=\Oh\Big(\frac{n^9}{\sigma^2}\Big). 
$$ 
\end{theorem} 

\proof
See Section~\ref{sec:main_proof}. 
\eproof

\section{Some properties of the condition number $\mu$}
\label{sec:condition_property}

There is a general geometric framework for defining condition numbers,
see~\cite[\S14.3]{Condition}. In our situation, 
this framework takes the following form. 

If $(A,\la,v)\in\V$ is well-posed, then the projection 
$\pi_1\colon\V\to\Cnn$ (cf.~(2)), around $(A,\la,v)$, 
has a local inverse $U\to\V, A\mapsto (A,G(A))$, that is defined on an open neighborhood $U$ of $A$ in $\Cnn$.
We call $G$ the {\em solution map}. The map $G$ decomposes as 
$G=(G_{\la},G_v)$ where 
$$
  G_{\la}:\Cnn\to\C
  \qquad{\mbox{and}}\qquad
  G_v:\Cnn\to\P(\C^n)
$$
associate to matrices $B$ in $U$ an eigenvalue and corresponding 
eigenvector, respectively. Let 
$$
  DG_{\la}:\Cnn\to\C
  \qquad{\mbox{and}}\qquad
  DG_v:\Cnn\to T_v
$$
be the derivatives of these maps at $A$ (which are linear maps 
between tangent spaces). The condition numbers for the eigenvalue 
$\la$ and the eigenvector $v$ of $A$ are defined as follows:
$$
  \mu_{\la}(A,\la,v):=\|DG_{\la}\|
  \qquad{\mbox{and}}\qquad
  \mu_v(A,\la,v):=\|DG_v\|
$$
where the norms are the operator norms with respect to the  
chosen norms (on $\Cnn$ we use the Frobenius norm and on $T_v$ 
the norm given by~\eqref{eq:defTip}).  
The following result, Lemma~14.17 in~\cite{Condition}, 
gives explicit descriptions of $DG_{\la}$ and 
$DG_v$. Before stating it, we recall that if $\la$ is an 
eigenvalue of $A$ there exists 
$u\in\C^n$ such $(A^*-\bar{\la}\Id)u=0$. We say that 
$u$ is a {\em left eigenvector} of $A$. 
Recall the linear map 
$A_{\lambda,v}\colon T_v\to T_v$ 
introduced in~\eqref{eq:defAlv}. 

\begin{lemma}\label{lem:14.17}
Let $(A,\la,v)\in\V$ and let $u$ be a left eigenvector of $A$ 
with eigenvalue~$\la$. Then:

{\bf (a)\ } We have $\langle v,u\rangle \ne 0$. 
\smallskip

{\bf (b)\ } If $\lambda$ is a simple eigenvalue of 
$A\in\C^{n\times n}$ with right eigenvector~$v$ 
and left eigenvector $u$, then the derivative of the solution 
map is given by 
$DG(A)(\dot{A}) = (\dot{\lambda},\dot{v})$, where 
\begin{equation}\tag*{\qed}
 \dot{\lambda} = \frac{\langle \dot{A}v,u\rangle}
  {\langle v,u\rangle}, \quad 
 \dot{v} = A_{\lambda,v}^{-1} \, P_{v^\perp} \dot{A} v .
\end{equation}
\end{lemma}

Lemma~\ref{lem:14.17} can be used to bound 
eigenvalue and eigenvector condition numbers.  
The following result is essentially Prop. 14.15 
in~\cite{Condition} (the only difference being that here 
we use Frobenius norms).

\begin{proposition}\label{pro:cn-eigenv}
Choosing the Frobenius norm on $T_A\Cnn = \C^{n\times n}$ and 
$\frac1{\|v\|}\, \|~\|$ on $T_v$, the condition numbers $\mu_v$ 
for the eigenvector problem and $\mu_{\la}$ for the eigenvalue 
problem satisfy: 
\begin{eqnarray*}
\mu_{\la}(A,\lambda,v)\;=\;\|DG_{\la}(A)\| 
& = & \frac{\|u\|\|v\|}{|\langle u ,v \rangle | },\\
\mu_{v}(A,\la,v)\;=\;\|DG_v(A)\| 
 & = & \big\| A_{\lambda,v}^{-1} \big\|. 
\end{eqnarray*}
\end{proposition}

\proof
For all $\dot{A} \in \C^{n\times n}$ we have,  
by the Cauchy-Schwarz inequality, 
$$
 |\langle \dot{A}v, u\rangle | \le \| \dot{A}v \| \|u\| 
 \le \|\dot{A}\|_F\|v\|\|u\|.
$$
This implies with Lemma~\ref{lem:14.17} that 
$$
 \|DG_{\la}(A)\| = \max_{\|\dot{A} \|_F= 1} 
 \frac{ |\langle \dot{A} v, u \rangle|}{|\langle v ,u \rangle|} 
 \leq \frac{\|v\|\|u\|}{|\langle v,u\rangle| } .
$$
Moreover, there exists  a rank one matrix $\dot{A}$ with 
$\|\dot{A}\| = \|\dot{A}\|_F = 1$ such that 
$\dot{A}v/\|v\| = u/\|u\|$, cf.\ \cite[Lemma 1.2]{Condition}. 
For this choice of $\dot{A}$ we have equality above. 
This proves the first assertion. 

In order to bound $\|DG_v(A)\|$ we note that for all $\dot{A}$,  
$$
 \big\| P_{v^\perp}\dot{A} v\big\| \ \le\  
 \big\|P_{v^\perp}\big\|\|\dot{A}\|\|v\| 
  \ \le\  \|\dot{A}\|_F\|v\|.
$$
Therefore, 
\begin{equation*}
 \max_{\|\dot{A} \|_F= 1} \big\| A_{\lambda,v}^{-1} 
 P_{v^\perp} \dot{A} v \big\| 
 \ \le\   \big\| A_{\lambda,v}^{-1} \big\| 
 \max_{\|\dot{A} \|_F= 1} \big\| P_{v^\perp}\dot{A} v\big\| \\
 \ \le\  \big\|A_{\lambda,v}^{-1}  \big\|\;\|v\|.
\end{equation*}
The inequality in the second assertion follows with 
Lemma~\ref{lem:14.17} (and the choice of norm $\|v\|^{-1} \|~\|$ 
on $T_v$). 
For the equality, let $w\in T_v$ be such that $\|w\|=\|v\|$. 
Again, by \cite[Lemma 1.2]{Condition}, there exists $\dot{A}$ such that 
$\|\dot{A}\|_F=1$ and $\dot{A} v= w$, hence 
$P_{v^\perp} \dot{A} v = w$. 
This implies
$$
  \max_{\|\dot{A}\|_F = 1} \| A_{\lambda,v}^{-1} \, P_{v^\perp} \dot{A} v \| 
  \ge \max_{w\in T_v \atop \|w\| = \|v\| } \| A_{\lambda,v}^{-1} \, w\|
 = \|v\| \,  \| A_{\lambda,v}^{-1} \| ,
$$
which completes the proof. 
\eproof

An immediate consequence of 
Proposition~\ref{pro:cn-eigenv} is that 
$\mu_v(A,\la,v)\; =\;\frac{\mu(A,\la,v)}{\|A\|_F}$. 
We next show that $\mu_{\la}(A,\la,v)$ can be similarly bounded 
in terms of $\mu(A,\la,v)$.

\begin{lemma}\label{le:ev_mu}
If $\la$ is a simple eigenvalue of $A$ with left eigenvector~$v$ and 
right eigenvector~$u$, then
$$
 \mu_{\la}(A,\la,v) 
\ \le\ \sqrt{1+\mu(A,\la,v)^2} .
$$
\end{lemma}

\proof
By unitary invariance, we may assume without loss of generality that 
$A=\begin{bmatrix}  \lambda & a \\ 0 & B\end{bmatrix}$ and $v=(1,0)$, 
where $a\in\C^{n-1}$ and $B\in\C^{(n-1)\times (n-1)}$. 
Clearly, $\|a\|\leq\|A\|_F$. 

We define $x:= ((\la\Id-B)^*)^{-1} a^*$. 
Then $a^* + B^* x = \bar{\lambda}x$, which implies 
$$
 \begin{bmatrix}  \bar{\lambda} & 0 \\ a^* & B^*\end{bmatrix}
 \begin{bmatrix}  1 \\ x \end{bmatrix} = 
 \bar{\lambda}\begin{bmatrix}  1 \\ x \end{bmatrix} .
$$
Hence $u:=(1,x)$ is a corresponding right eigenvector 
of $\la$. Also, 
$$
    \|x\| \le  \|((\la\Id -B)^*)^{-1}\| \cdot \|a^*\| 
           = \|(\la\Id -B)^{-1}\| \cdot \|a\|  \le \mu(A,\la,v)
$$
since $A_{\la,v}= B-\la\Id$ and $\|a\| \le \|A\|_F$. 
Now, 
$$
 \frac{\|u\| \|v\|}{|\langle u,v\rangle|} = \|u\| = \sqrt{1 + \|x\|^2}  
$$
and the assertion follows with Proposition~\ref{pro:cn-eigenv}. 
\eproof
\medskip

\proofof{Proposition~\ref{prop:dotzdotf-spherical}}
The first two inequalities are immediate from 
Lemma~\ref{le:ev_mu} and 
$\mu_v(A,\la,v) = \frac{\mu(A,\la,v)}{\|A\|_F}$. 
For the third we have 
\begin{equation*}
\|\dot{\Gamma}(t)\|\;=\;\|(\dot{A},\dot{\la},\dot{v})\| 
 \;\le\;\|\dot{A}\|\sqrt{1+\mu_t^2+(1+\mu_t^2)} 
  \;\leq\; \|\dot{A}\|\sqrt{6\mu_t^2} 
\end{equation*}
the last inequality since $\mu_t\geq\frac{1}{\sqrt{2}}$ 
(Lemma~\ref{le:lb_mu}).
\eproof

\section{The Lipschitz property for the eigenpair problem}
\label{se:lipschitz} 

In this section we prove Theorem~\ref{th:lipsch}. Recall, 
we assume that $(A,\la,v)$ is in the solution
variety, i.e., $Av=\lambda v$, but we {\em do not} 
require that $A'v'=\la'v'$.
The following result is the main stepping stone.

\begin{proposition}\label{prop:lipsch1}
Let $A,A'\in\C^{n\times n}$ be such that $\|A\|_F = \|A'\|_F = 1$, let 
$v,v'\in\C^n$ be nonzero, and let  $\la,\la'\in\C$ be such that 
$Av=\lambda v$.  
Suppose that
$$
 \mu(A,\la,v) \big( \|A'-A\|_F + |\la' - \la| + \dpr(v',v) \big) \ \le\  
 \frac{\e}{7.2}
$$
for $0<\e\leq 0.37$. Then we have 
$$
 \frac{1}{1+\e}\,\mu (A,\lambda,v)\ \le\ \mu(A',\la',v') \ \le\  
 (1+\e)\, \mu(A,\la,v) .
$$
\end{proposition}
 
Before proceeding, we note that Theorem~\ref{th:lipsch} is an 
immediate consequence of this proposition, since 
\begin{equation*} 
\dist((A,\la,v),(A',\la',v'))\leq \frac{1}{\sqrt{3}}
\Big(\|A'-A\|_F + |\la' - \la| + \dpr(v',v)\Big).
\end{equation*}

We shall provide the proof of Proposition~\ref{prop:lipsch1}
in several steps. We begin with the following easy observation 
whose proof is left to the reader 
(cf.~\cite[Lemma~16.40]{Condition}). 

\begin{lemma}\label{le:OP}
Suppose that $v',v\in\C^n$ are nonzero with $\d:=\dpr(v',v) < \pi/2$. 
Then the restriction 
$P_{v^\perp} |_{T_{v'}}\colon T_{v'} \to T_v$ 
of $P_{v^\perp}$ to $T_{v'}$ is invertible and we have 
\begin{equation}\tag*{\qed}
 \| P_{v^\perp} |_{T_{v'}}\| = \cos\d, \quad 
  \| (P_{v^\perp} |_{T_{v'}})^{-1}\| = (\cos\d)^{-1}.
\end{equation}
\end{lemma}
\medskip

In a first step, we fix $v$ and only perturb $A$ and $\la$. 

\begin{lemma}\label{le:uno}
Suppose that $A,A'\in\C^{n\times n}$ are such that 
$\|A\|_F = \|A'\|_F = 1$, $\la,\la'\in\C$ and $v\in\C^n$ is nonzero.
Then,  
$$
 \mu(A,\la,v) \big( \|A'-A\|_F + |\la' - \la| \big) \ \le\ \e < 1 
$$
implies that 
$$
(1-\e)\,\mu (A,\lambda,v)\ \le\ \mu(A',\la',v) \ \le\  
 \frac{1}{1-\e}\, \mu(A,\la,v) .
$$
\end{lemma}

\proof
Recall that $A_{\la,v} =P_{v^\perp}(A-\la\Id)$, seen as a linear
endomorphism of~$T_v$. 
We have $A'_{\la',v} = A_{\la,v} + \Delta$, where 
$\Delta:= P_{v^\perp} ((A-A')+(\la'-\la)\Id)$, interpreted as a linear
endomorphism of~$T_v$. Note that 
$$
 \|\Delta\|\;\le\;\|P_{v^\perp}\| \cdot 
  \big(\|A-A'\|+\|(\la' -\la)\Id\|\big) 
 \;\le\;\|A-A'\|+|\la' -\la|\;\le\;\|A-A'\|_F+|\la' -\la|.
$$
Our assumption implies 
$\|(A_{\la,v})^{-1}\|\cdot \|\Delta\| \ \le \e < 1$ 
(note that 
$\mu(A,\la,v) = \| (A_{\la,v})^{-1}\|$ 
since we assume that $\|A\|_F=1$).  
Also, Lemma~15.7 in~\cite{Condition} implies that 
$$
 \| (A'_{\la',v})^{-1}\| \ \le\ \| (A_{\la,v})^{-1}\| \cdot 
  \frac{1}{1 - \| (A_{\la,v})^{-1}\|\cdot \|\Delta \|} .
$$
Using again $\|A\|_F=\|A'\|_F=1$ it follows that 
\begin{equation}\label{eq:UG1}
\mu(A',\la',v) \ \le\  \frac{1}{1-\e}\, \mu(A,\la,v) .
\end{equation}
It now suffices to prove that 
$$
 (1-\e)\,\mu (A,\lambda,v)\ \le\ \mu(A',\la',v).
$$
Since this is trivial in the case $\mu(A,\la,v) \le \mu(A',\la',v)$, we may 
assume that 
$\mu(A',\la',v) <\mu(A,\la,v)$. Then, by our assumption, 
$$
 \mu(A',\la',v) \big( \|A'-A\|_F + |\la' - \la| \big) \ \le\  
 \mu(A,\la,v) \big( \|A'-A\|_F + |\la' - \la| \big) \le  \e < 1 .
$$
Hence we can apply \eqref{eq:UG1} by switching the roles of $(A',\la',v)$ and $(A,\la,v)$. 
This gives 
$$
 \mu(A,\la,v) \ \le\  \frac{1}{1-\e}\, \mu(A',\la',v) , 
$$
which completes the proof.
\eproof

In a second step, we fix $A$ and $\lambda$ and only perturb $v$. 
Now we need to assume one of the triples in $\V$.

\begin{lemma}\label{le:due}
Let $A\in\C^{n\times n}$ with $\|A\|_F = 1$, $\la\in\C^n$ and 
$v,v'\in\C$ be nonzero such that $Av = \la v$. 
Then
$$
  3.6\,\mu(A,\la,v)\, \dpr(v,v') \ \le\  \e < 1 
$$
implies that 
$$
\frac{\mu(A,\la,v)}{1+\e} \ \le\ \mu (A,\lambda,v')\ \le\ \frac{\mu(A,\la,v)}{1-\e} .
$$
\end{lemma}

\proof
Since $\mu(A,\la,v) < \infty$, $\la$ is a simple eigenvalue of $A$. 
Let $u$ denote a corresponding right eigenvector of $A$, that is, 
$A^*u = \bar{\lambda}u$. We have 
$$
 \langle u, (\la\Id -A) x\rangle = \langle u , \la x \rangle - \langle u,  
  Ax \rangle 
 = \bar{\la} \langle u ,x \rangle - \langle A^* u,  x \rangle = 0 ,
$$ 
hence the image of $\la\Id -A $ is contained in $T_u$ and thus equals 
$T_u$ for reasons of dimension. 

Let $\pi_v\colon T_u \to T_v$ denote the restriction to $T_u$ of the
orthogonal projection~$P_{v^\perp}$.  Since $\langle u, v \rangle\ne
0$, the map $\pi_v$ is an isomorphism. We denote by $\pi_{v'}\colon
T_u \to T_{v'}$ the restriction to $T_u$ of the orthogonal
projection~$P_{v'^{\perp}}$ and we define $\gamma := \pi_{v'}\circ
\pi_v^{-1}$.  Moreover, we let $\pi\colon T_{v'}\to T_v$ denote the
restriction to $T_{v'}$ of the orthogonal projection~$P_{v}$.  We
further write $\Phi:=\la\Id-A$ and consider the following commutative
diagram: 
\begin{center}
\begin{tikzpicture}
\path (-0.25,0.3) node[right]{$T_u$};
\path (1.5,0.3) node[right]{$T_v$};
\path (-1.9,0.3) node[right]{$T_v$};
\draw[->] (0.15,0) -- (1.45,-1) node[above =3mm]{};
\path (0.8,-0.5) node[right]{\small $\pi_{v'}$};
\draw[->] (1.9,0) -- (1.9,-1) node[above =3mm]{};
\path (1.85,-0.45) node[right]{\small $\gamma$};
\draw[->] (-1.6,-1) -- (-1.6,0) node[above =3mm]{};
\path (-2.05,-0.5) node[right]{\small $\pi$};
\draw[->] (-1.2,0.3) -- (-0.3,0.3) node[above =3mm]{};
\path (-1.,0.5) node[right]{\small $\Phi$};
\draw[->] (0.4,0.3) -- (1.5,0.3) node[above =3mm]{};
\path (0.7,0.5) node[right]{\small $\pi_v$};
\path (-1.9,-1.3) node[right]{$T_{v'}$};
\draw[->] (-1.2,-1.3) -- (1.45,-1.3) node{};
\path (-0.45,-1) node[right]{\small $A_{\la,v'}$};
\path (1.5,-1.3) node[right]{$T_{v'}$.};
\end{tikzpicture}
\end{center}
%
The top map $\pi_v \circ \Phi \colon T_v\to T_v$ equals $A_{\la,v}$ by
definition.  Our assumption $Av=\la v$ means $\Phi(v)=0$. Hence
$\Phi(x) = \Phi(\pi(x))$ for $x\in\C^n$.  This implies that the bottom
map $\pi_{v'}\circ\Phi\circ\pi$ indeed equals $A_{\la,v'}$.  We
conclude that 
$$ 
  A_{\la,v'} = \gamma\circ A_{\la,v} \circ\pi,\quad
  A_{\la,v'}^{-1} = \pi^{-1}\circ A_{\la,v}^{-1}\circ\gamma^{-1},\quad
  \|A_{\la,v'}^{-1}\| \le \|\pi^{-1}\| \cdot \|\gamma^{-1}\|
  \cdot\|A_{\la,v}^{-1}\| .  
$$ 
We will see in a moment that $\pi_{v'}^{-1}$ is bijective. 
Then, $\gamma^{-1} =\pi_v\circ\pi_{v'}^{-1}$, 
hence $\|\gamma^{-1}\| \le \|\pi_v\|\cdot \|\pi_{v'}^{-1}\|$.

We use the abbreviations $\d:=\dpr(v',v)$, $\d_0:=\dpr(u,v)$,  and 
$\d_1:=\dpr(u,v')$.
Note that $\d_0-\d_1 \le \d_1\le \d_0 + \d$ by the triangle inequality. 
Let us proceed with some estimates. 

Using the bound $\cos\d\ge 1- \frac{2}{\pi}\d$, which is valid for 
$\d\le\pi/2$, we get 
\begin{equation}\label{eq:cosAT}
 \frac{\cos(\d_0 + \d)}{\cos\d_0} = \cos\d - \tan(\d_0) \sin(\d) \ \ge\  
 1 - \Big(\frac{2}{\pi} + \tan\d_0 \Big) \d .
\end{equation}
This implies
$$
 \cos\d \cdot \frac{\cos(\d_0 + \d)}{\cos\d_0} \ \ge\ 
 1 - \Big(\frac{4}{\pi} + \tan\d_0 \Big) \d .
$$
We write $\mu:=\mu(A,\la,v)$ to simplify the notation.  
Lemma~\ref{le:ev_mu} provides the following estimate
\begin{equation}\label{eq:tanES}
 \tan\d_0 \ \le\  \frac{1}{\cos\d_0} = \sqrt{1+\mu^2} \ \le\ \sqrt{3}\, 
\mu ,
\end{equation}
where the last inequality is due to Lemma~\ref{le:lb_mu}. 
Again using Lemma~\ref{le:lb_mu}, we estimate 
$$
 \frac{4}{\pi} + \tan\d_0 \ \le\  
 \Big(\frac{4\sqrt{2}}{\pi} +\sqrt{3} \Big) \, \mu\ \le\  3.6\, \mu .
$$
We conclude that 
\begin{equation}\label{eq:ESM}
 \cos\d \cdot \frac{\cos(\d_0 + \d)}{\cos\d_0} \ \ge\ 
 1 - 3.6\, \mu \d >0 ,
\end{equation}
where the positivity is a consequence of our assumption. 
This shows that $\d_1 \le \d_0 + \d < \pi/2$ and hence 
$\pi_{v'}$ is indeed bijective. 

Lemma~\ref{le:OP} yields 
$$
 \|\pi_v\| \le \cos\d_0, \quad 
 \|\pi^{-1}\| \le \frac{1}{\cos\d}, \quad 
 \|\pi_{v'}^{-1}\| \le \frac{1}{\cos\d_1} \le  \frac{1}{\cos(\d_0 + \d)} .
$$
Recall that 
$\|A_{\la,v'}^{-1}\| \le \|\pi^{-1}\| \cdot \|\pi_v\|\cdot \|\pi_{v'}^{-1}\| 
\cdot\|A_{\la,v}^{-1}\|$.
We obtain 
$$
 \|A_{\la,v'}^{-1}\| \ \le\ 
 \frac{1}{\cos\d}\cdot \frac{\cos\d_0}{\cos(\d_0 + \d)} 
 \cdot \|A_{\la,v}^{-1}\| .
$$
The estimate \eqref {eq:ESM} yields 
\begin{equation}\label{eq:1st-ES}
 \|A_{\la,v'}^{-1}\| \ \le\ 
 \frac{1}{1 - 3.6\mu \d} \cdot \|A_{\la,v}^{-1}\| .
\end{equation}

For the other inequality, we proceed similarly. We have 
$$
   A_{\la,v} = \gamma^{-1}\circ A_{\la,v'} \circ\pi^{-1},\quad 
  A_{\la,v}^{-1} = \pi\circ A_{\la,v'}^{-1}\circ\gamma,\quad 
  \|A_{\la,v}^{-1}\| \le \|\pi\| \cdot \|\gamma\| \cdot\|A_{\la,v'}^{-1}\| .
$$
Moreover, $\|\gamma\| \le \|\pi_{v'}\|\cdot \| \pi_{v}^{-1}\|$. 
Lemma~\ref{le:OP} yields 
$$
 \|\pi_{v'}\| \le \cos\d_1 \le \cos(\d_0-\d), \quad 
 \|\pi_v^{-1}\| \le \frac{1}{\cos\d_0}, \quad 
 \|\pi\| \le \cos\d .
$$
This implies
$$
 \|A_{\la,v}^{-1}\| \ \le\ \cos\d\cdot \frac{\cos(\d_0-\d)}{\cos\d_0} 
 \cdot \|A_{\la,v'}^{-1}\| .
$$ 
Using~\eqref{eq:cosAT}  and \eqref{eq:tanES} we estimate
$$
 \frac{\cos(\d_0-\d)}{\cos\d_0} = \cos\d' + \tan(\d_0) \sin(\d) \ \le\ 
  1 + \sqrt{3}\, \mu\, \d . 
$$
We arrive at
$$
\|A_{\la,v}^{-1}\| \ \le\ (1 + \sqrt{3}\, \mu\, \cdot\d)  \|A_{\la,v'}^{-1}\| ,
$$
which, together with~\eqref{eq:1st-ES} and $\|A\|_F=1$, 
completes the proof.
\eproof
\medskip
 
\proofof{Proposition~\ref{prop:lipsch1}}
The assumption 
\begin{equation}\label{eq:ass}
 \mu(A,\la,v) \big( \|A'-A\|_F + |\la' - \la| + \dpr(v',v) \big) 
 \ \le\  \frac{\e}{7.2}
\end{equation}
yields $3.6\, \mu(A,\la,v) \dpr(v',v) \le \frac{\e}{2}$. 
Lemma~\ref{le:due} implies that 
\begin{equation}\label{eq:ass0}
 \mu(A,\la,v') \ \le\ \frac{\mu(A,\la,v)}{1-\e/2} 
 \ \le\ 2\, \mu(A,\la,v) , 
\end{equation}
where we have used $\e\le 1$ for the last equality. 
Combining this with~\eqref{eq:ass}, we conclude that 
\begin{equation}\label{eq:ass2}
 \mu(A,\la,v')  \big( \|A'-A\|_F + |\la' - \la| \big)  \ \le\ \frac{\e}{3.6}. 
\end{equation}
Lemma~\ref{le:uno} and inequality~\eqref{eq:ass0} imply that 
$$
 \mu(A',\la',v') \ \le\ \frac{\mu(A,\la,v')}{1-\e/3.6} \ \le\ 
 \frac{\mu(A,\la,v)}{(1-\e/3.6)(1-\e/2)} \ \le\ 
 (1+\e)\mu(A,\la,v),
$$
the last since $\e\leq 0.37$. This 
proves one of the claimed inequalities.

For the other inequality we proceed similarly. We use again 
Lemma~\ref{le:due} now to obtain  
\begin{equation}\label{eq:ass33}
  \frac{\mu(A,\la,v)}{1+\e/3.6} \ \le\ \mu(A,\la,v') .
\end{equation}
Inequality~\eqref{eq:ass2} allows us to use again 
Lemma~\ref{le:uno}, now to obtain
\begin{equation}\label{eq:lip81}
 (1- \e/2)\, \mu(A,\la,v') \ \le\ \mu(A',\la',v') .
\end{equation}
Combining~\eqref{eq:ass33} and~\eqref{eq:lip81} we 
obtain
$$
  \mu(A',\la',v') \geq (1- \e/2)\, \mu(A,\la,v') 
  \geq \frac{1- \e/2}{1+\e/3.6} \, \mu(A,\la,v)
  \geq \frac{1}{1+\e}\, \mu(A,\la,v)
$$
the last inequality, again, since $\e\leq 0.37$. 
This completes the proof. 
\eproof

\section{Proof of Proposition~\ref{thm:main_path_following}}
\label{sec:homotopy}

The proof relies on three main ingredients: 
Proposition~\ref{prop:dotzdotf-spherical}, and 
Theorems~\ref{th:lipsch} and~\ref{th15.1}. Different versions 
of it have been written down for different contexts: 
complex polynomial systems (\cite[Theorem~3.1]{BuCu11} 
or~\cite[Theorem~17.3]{Condition}), 
ditto but with finite precision~(\cite[Theorem~4.3]{BriCuPeRo}), 
and Hermitian matrices~(\cite[Proposition~4]{ArmCuc}).  

We set $\e:=0.12$ and $C_\e:=\frac{\e}{12.5}=0.0096$.  
Furthermore, we let  
$\xi:=\frac{2C_\e(1-\e)}{3\sqrt{6}(1+\e)^4}\approx 0.001461$. 
Also, as in \S\ref{subsec:ALH}, we write $\a=\dS(M,A)$. 

Let $0=\tau_0<\tau_1<\cdots<\tau_K=1$ and
$(\lambda_0,v_0)=(\zeta_0,w_0),(\zeta_1,w_1),\ldots,(\zeta_K,w_K)$
be the sequences of $\tau$-values and pairs in $\C\times\C^n$
generated by the algorithm \EC. We simplify notation and 
write $Q_i$ instead of $Q_{\tau_i}$ and
$(\lambda_i,v_i)$ instead of $(\lambda_{\tau_i},v_{\tau_i})$. 
(There is no danger of confusing this $(\lambda_1,v_1)$ with the one 
appearing in the statement of Proposition~\ref{thm:main_path_following}.)
We associate with the solution path~\eqref{eq:curva} in~$\V$ 
the following curve in $\IS(\Cnn)\times\C\times\C^n$: 
\begin{equation}\label{eq:curva2}
  [0,1]\to V,\quad \tau \mapsto   (P_\tau,\widehat{\lambda_\tau},v_\tau) 
   := \Big(\frac{Q_\tau}{\|Q_\tau\|_F},
       \frac{\lambda_\tau}{\|Q_\tau\|_F},v_\tau\Big).
\end{equation}
We also write $P_i$ instead of $P_{\tau_i}$. 
The meaning of the parameterization by~$\tau$ is that 
$\a\tau$ is the parameterization of $\tau\mapsto P_\tau$ by arc length, 
which means that 
$\Big\|\frac{d P_\tau}{d\tau}\Big\| = \a$. 

We will carry out the proof on the curve~\eqref{eq:curva2}  
in the sphere $\IS(\Cnn)$. We do so to simplify the exposition 
and without implying that algorithm~\EC\ should be 
modified to normalize matrices. Indeed, all the functions 
on triples in $\Cnn\times\C\times\C^n$ involved in our proof 
---$\dist$, $\mu$, and Newton's operator--- 
are scale invariant on the first two components. Furthermore, to avoid 
burdening the notation, we will write $\lambda$ instead of 
$\widehat{\lambda}$. This should introduce no confusion. 

The following result is the technical core in the proof 
of Proposition~\ref{thm:main_path_following}. 

\begin{proposition}\label{prop:ALH_induction}
For $i=0,\ldots,K-1$, the following statements
are true:

\noindent {\bf (a)}\quad $\displaystyle
 \dist((P_i,\zeta_i,w_i),(P_i,\lambda_i,v_i))
 \leq \frac{C_\e}{\mu(P_i,\lambda_i,v_i)}$.
\smallskip

\noindent {\bf (b)}\quad $\displaystyle
    \frac{\mu(P_i,\zeta_i,w_i)}{1+\e}\leq
    \mu(P_i,\lambda_i,v_i)\leq (1+\e)\mu(P_i,\zeta_i,w_i)$.
\smallskip

\noindent {\bf (c)}\quad $\displaystyle
 \dist((P_i,\lambda_i,v_i),(P_{i+1},\lambda_{i+1},v_{i+1}))
 \leq \frac{C_\e}{\mu(P_i,\lambda_i,v_i)}\frac{2(1-\e)}{3(1+\e)}$.
\smallskip

\noindent {\bf (d)}\quad $\displaystyle
 \dist((P_{i+1},\zeta_i,w_i),(P_{i+1},\lambda_{i+1},v_{i+1})) \leq 
 \frac{2C_\e}{(1+\e)\mu(P_i,\lambda_i,v_i)}$. 
\smallskip

\noindent {\bf (e)}\quad $(\zeta_i,w_i)$ is an approximate eigenpair 
of $P_{i+1}$ with associated eigenpair $(\lambda_{i+1},v_{i+1})$.
\eproof
\end{proposition}

\proof
We proceed by induction, showing that
$$
  ({\bf a},i)\Rightarrow ({\bf b},i)\Rightarrow
  ({\bf c},i) 
  \Rightarrow({\bf d},i)\Rightarrow \big( ({\bf e},i) \mbox{\rm\ and } 
  ({\bf a},i+1) \big).
$$
Inequality (a) for $i=0$ is trivial.

Assume now that (a) holds for some $i \le K-1$. Then,
Theorem~\ref{th:lipsch} (with $A=A'=P_i$) implies
$$
   \frac{\mu(P_i,\zeta_i,w_i)}{1+\e}\leq
   \mu(P_i,\lambda_i,v_i)\leq (1+\e)\mu(P_i,\zeta_i,w_i)
$$
and thus (b). We now prove (c). 
To do so, let $\tau_*>\tau_i$ be such that
$$
 \int_{\tau_i}^{\tau_*}\bigg\|\frac{d(P_\tau,\lambda_\tau,v_\tau)}
    {d\tau}\bigg\| d\tau
 =\frac{C_\e}{\mu(P_i,\lambda_i,v_i)}\,\frac{2(1-\e)}{3(1+\e)}
$$ 
or $\tau_*=1$, whichever is smaller.
Then, for all $t\in[\tau_i,\tau_*]$,
\begin{equation}\label{eq:EC_uno}
\begin{split}
 \dist((P_i,\lambda_i,v_i),(P_t,\lambda_t,v_t)) 
&\leq \dR((P_i,\lambda_i,v_i),(P_t,\lambda_t,v_t))\\[3pt] 
&= \int_{\tau_i}^t\bigg\|\frac{d(P_\tau,\lambda_\tau,v_\tau)}
    {d\tau}\bigg\|\, d\tau    
    \leq \frac{C_\e}{\mu(P_i,\lambda_i,v_i)}\,\frac{2(1-\e)}{3(1+\e)}, 
\end{split}
\end{equation}
the first inequality by~\eqref{eq:chord}. 
It is therefore enough to show that $\tau_{i+1}\leq \tau_*$. This is
trivial if $\tau_*=1$. We therefore assume $\tau_*<1$. The 
bound above allows us to apply Theorem~\ref{th:lipsch} and
to deduce, for all $\tau\in[\tau_i,\tau_*]$,
\begin{equation}\label{eq:step_bound}
   \frac1{1+\e}\,\mu(P_i,\lambda_i,v_i) \leq
    \mu(P_\tau,\lambda_\tau,v_\tau) \leq (1+\e)\,\mu(P_i,\lambda_i,v_i).
\end{equation}
Proposition~\ref{prop:dotzdotf-spherical} implies that 
$$
 \bigg\|\frac{d}{d\tau} (P_\tau,\lambda_\tau,v_\tau)\bigg\| 
 \le \sqrt{6}\, \mu(P_\tau,\lambda_\tau,v_\tau)\, 
 \bigg\|\frac{d}{d\tau} P_\tau\bigg\| 
$$
We now use~\eqref{eq:step_bound} to deduce that
\begin{equation*}
\begin{split}
 \frac{C_\e}{\mu(P_i,\lambda_i,v_i)} \frac{2(1-\e)}{3(1+\e)} 
  \ =\  \int_{\tau_i}^{\tau_*} \bigg\|\frac{d(P_\tau,\lambda_\tau,v_\tau)}
    {d\tau}\bigg\| d\tau
  \ \leq\ \int_{\tau_i}^{\tau_*}\sqrt{6}\,\mu(P_\tau,\lambda_\tau,v_\tau)
\bigg\|\frac{d}{d\tau} P_\tau\bigg\| d\tau \\
  \leq\  \sqrt{6}\, (1+\e)\, \mu(P_i,\lambda_i,v_i) 
  \int_{\tau_i}^{\tau_*} \bigg\|\frac{d}{d\tau} P_\tau\bigg\| d\tau 
  \ = \  \sqrt{6}\, (1+\e)\, \mu(P_i,\lambda_i,v_i)\, \dS(P_i,P_{\tau_*}) . 
\end{split}
\end{equation*}
Consequently, using (b), we obtain
$$
  \dS(P_i,P_{\tau_*}) \geq \frac{2C_\e(1-\e)}{3\sqrt{6}(1+\e)^2
  \mu^2(P_i,\lambda_i,v_i)}
  \geq \frac{2C_\e(1-\e)}{3\sqrt{6}(1+\e)^4\mu^2(P_i,\zeta_i,w_i)} .
$$
Recall that the parameter $\xi$ in \EC\ was chosen as 
$\xi = \frac{2C_\e(1-\e)}{3\sqrt{6}(1+\e)^4}$.
By the definition of $\tau_{i+1}-\tau_i$ in \EC\ we have
$\a(\tau_{i+1}-\tau_i) = \frac{\xi}{\mu^2(P_i,\zeta_i,w_i)}$.
So we obtain
$$
 \dS(P_i,P_{\tau_*}) \geq \a(\tau_{i+1}-\tau_i)\,=\,\dS(P_i,P_{i+1}).
$$
This implies $\tau_{i+1}\leq \tau_*$ as claimed, and hence 
inequality (c) follows from \eqref{eq:EC_uno} with $t=\tau_{i+1}$. 
With it, we may apply
Theorem~\ref{th:lipsch} once more to deduce, for all
$\tau\in[\tau_i,\tau_{i+1}]$,
\begin{equation}\label{eq:b1}
    \frac{\mu(P_i,\lambda_i,v_i)}{1+\e} \leq \mu(P_\tau,\lambda_\tau,v_\tau)
    \leq (1+\e) \mu(P_i,\lambda_i,v_i).
\end{equation}
We now observe that 
\begin{eqnarray*}
   \dist((P_{i+1},\zeta_i,w_i),(P_i,\zeta_i,w_i)) &=& \|P_{i+1}-P_i\|_F 
   \;\leq\; 
   \dS(P_i,P_{i+1})\;=\;\a(\tau_{i+1}-\tau_i) \\
  &=&\frac{\xi}{\mu^2(P_i,\zeta_i,w_i)} 
  \;\leq\; \frac{2C_\e(1-\e)}{3\sqrt{6}(1+\e)\mu(P_i,\zeta_i,w_i)}
\end{eqnarray*}
and use this bound, together with the triangle 
inequality,~\eqref{eq:step_bound},~(a), and (c) 
to obtain
\begin{equation}\label{eq:EC_tre}
\begin{split}
  \dist((P_{i+1},&\zeta_i,w_i),(P_{i+1},\lambda_{i+1},v_{i+1})) 
\leq \dist((P_{i+1},\zeta_i,w_i),(P_i,\zeta_i,w_i)) \\
&\qquad +\dist((P_i,\zeta_i,w_i),(P_i,\lambda_i,v_i)) \\
&\qquad +\dist((P_i,\lambda_i,v_i),(P_{i+1},\lambda_{i+1},v_{i+1})) \\
 &\leq\frac{2C_\e(1-\e)}{3\sqrt{6}(1+\e)\mu(P_i,\zeta_i,w_i)} 
   +\frac{C_\e}{\mu(P_i,\lambda_i,v_i)}
    +\frac{C_\e}{\mu(P_i,\lambda_i,v_i)}\frac{2}{3}\frac{1-\e}{1+\e} \\
  &\leq \frac{2C_\e}{(1+\e)\mu(P_i,\lambda_i,v_i)}
\end{split}
\end{equation} 
which proves (d). We now note that $\frac{2C_\e}{1+\e}<c_0 = 0.2881$. 
We can therefore apply Theorem~\ref{th15.1} 
to deduce that $(\zeta_i,w_i)$ is an approximate eigenpair
of~$P_{i+1}$ associated with its eigenpair $(\lambda_{i+1},v_{i+1})$, 
and hence (e) holds.

It follows from (e) that 
$(\zeta_{i+1},w_{i+1})=N_{P_{i+1}}(\zeta_i,w_i)$ satisfies
$$
\dist((P_{i+1},\zeta_{i+1},w_{i+1}),(P_{i+1},\lambda_{i+1},v_{i+1}))
\leq \frac12\,\dist((P_{i+1},\zeta_i,w_i),(P_{i+1},\lambda_{i+1},v_{i+1})).
$$
Using this bound, (d) and the right-hand inequality in~\eqref{eq:b1} 
with $\tau=\tau_{i+1}$, we obtain 
$$
   \dist((P_{i+1},\zeta_{i+1},w_{i+1}),(P_{i+1}\lambda_{i+1},v_{i+1}))
  \leq \frac{C_\e}{(1+\e)\mu(P_i,\lambda_i,v_i)}\\
  \leq \frac{C_\e}{\mu(P_{i+1},\lambda_{i+1},v_{i+1})} ,
$$
which proves (a) for $i+1$. The proposition is thus proved.
\eproof

The following is an immediate consequence of 
Proposition~\ref{prop:ALH_induction}(c) and 
Theorem~\ref{th:lipsch}. 

\begin{corollary}\label{cor:b1}
For all $i=0,\ldots,K-1$ and $\tau\in[\tau_i,\tau_{i+1}]$,
we have
\begin{equation}\tag*{\qed}
\frac{\mu(Q_i,\lambda_i,v_i)}{1+\e} \leq \mu(Q_\tau,\lambda_\tau,v_\tau)
    \leq (1+\e) \mu(Q_i,\lambda_i,v_i).
\end{equation}
\end{corollary}

\proofof{Proposition~\ref{thm:main_path_following}}
It follows from Proposition~\ref{prop:ALH_induction}(e) 
for $i=K-1$ that $(\zeta_{K-1},w_{K-1})$ is an 
approximate eigenpair of $Q_K=A$ with associated 
eigenpair~$(\lambda_K,v_K)$.  
Consequently, so is the 
returned point~$(\zeta_K,w_K)=N_{A}(\zeta_{K-1},w_{K-1})$.

Consider now any $i\in\{0,\ldots,K-1\}$. 
Using  Corollary~\ref{cor:b1}, 
and Proposition~\ref{prop:ALH_induction}(b), 
and by the choice of the step size $\Delta\tau$ 
in Algorithm~\ref{alg:EC}, we obtain
\begin{eqnarray*}
  \int_{\tau_i}^{\tau_{i+1}} \mu^2(Q_\tau,\lambda_\tau,v_\tau)d\tau
  &\geq&
  \int_{\tau_i}^{\tau_{i+1}
  }\frac{\mu^2(Q_i,\lambda_i,v_i)}{(1+\e)^2}d\tau
  \;=\;  \frac{\mu^2(Q_i,\lambda_i,v_i)}{(1+\e)^2} (\tau_{i+1}-\tau_i) \\
  &\geq&\frac{\mu^2(Q_i,\zeta_i,w_i)}{(1+\e)^4}(\tau_{i+1}-\tau_i)\\
  &=& \frac{\mu^2(Q_i,\zeta_i,w_i)}{(1+\e)^4}
       \frac{\xi}{\a \mu^2(Q_i,\zeta_i,w_i)}
  \;=\;  \frac{\xi}{(1+\e)^4\a } \\
   &\ge& \frac1{1077\,\a}.
\end{eqnarray*}
This implies 
$$
  \int_0^1 \mu^2(Q_\tau,\lambda_\tau,v_\tau)d\tau
   \ge \frac{K}{1077\,\a},
$$
which proves the stated upper bound on~$K$. 
The lower bound follows from 
\begin{align}
  \int_{\tau_i}^{\tau_{i+1}} \mu^2(Q_\tau,\lambda_\tau,v_\tau)d\tau
  \ \leq\ &
  \int_{\tau_i}^{\tau_{i+1}}
   \mu^2(Q_i,\lambda_i,v_i)(1+\e)^2d\tau
\notag \\
  =\ &  \mu^2(Q_i,\lambda_i,v_i)(1+\e)^2(\tau_{i+1}-\tau_i) 
\notag \\
  \leq\ &\mu^2(Q_i,\zeta_i,w_i)(1+\e)^4(\tau_{i+1}-\tau_i)
\notag \\
  =\ & \frac{\xi(1+\e)^4}{\a }
   \ \le\  \frac{1}{434\,\a}.\tag*{\qed}
\end{align}
As pointed out at the beginning of this section, the proof 
above follows the steps of the one in~\cite[Proposition~4]{ArmCuc}. 
In both cases we obtain 
$$
   \frac{3}{C_\e(1-\e)} \leq \frac{K}{\a\int_0^1 
   \mu^2(Q_\tau,\lambda_\tau,v_\tau)d\tau}
   \leq \frac{3(1+\e)^8}{C_\e(1-\e)}.  
$$
The difference in the actual constants in both statements 
is due to the difference between the values of $\e$ and $C_\e$ 
in (the inadecuate)~\cite[Proposition~3.22]{Armentano:13} 
used in~\cite{ArmCuc} 
and its improved version, Theorem~\ref{th:lipsch}  
(see the paragraph following Proposition~3 in~\cite{ArmCuc}). 

\section{Proof of Theorem~\ref{th:mu2-bound}}\label{se:muave} 

This section is the technical heart of the paper. We divide it into 
several subsections, the first of which summarizes some notions 
and tools of probability theory on Riemannian manifolds.

\subsection{The coarea formula}

On a Riemannian manifold~$M$ 
there is a well-defined measure $\vol_M$ obtained 
by integrating the indicator functions $\uno_A$ 
of Borel-measurable subsets $A\subseteq M$ 
against the volume form $dM$ of~$M$: 
$$
    \vol_M(A) = \int_M  \uno_A \, dM. 
$$  
Dividing $\uno$ by $\vol_M(M)$ if $\vol_M(M)<\infty$, 
this leads to a natural notion of {\em uniform distribution}
on $M$, which we will denote by $\msU(M)$. 
More generally, we will call any measurable function 
$f\colon M\to [0,\infty]$ such that $\int_M f\, dM = 1$
a {\em probability density} on~$M$.

The coarea formula is an extension of the transformation formula  
to not necessarily bijective smooth maps 
between Riemannian manifolds. 
In order to state it, we first need to generalize 
the notion of Jacobians.

Suppose that $M,N$ are Riemannian manifolds of 
dimensions~$m$, $n$,
respectively such that $m\ge n$. Let $\psi\colon M\to N$ be
a smooth map. By definition, the derivative $D\psi(x)\colon
T_xM\to T_{\psi(x)} N$ at a regular point $x\in M$ is surjective.
Hence the restriction of $D\psi(x)$ to the orthogonal complement of
its kernel yields a linear isomorphism. The absolute value of its
determinant is called the
{\em normal Jacobian} of $\psi$ at $x$ and 
denoted by $\NJ\psi(x)$.
We set $\NJ\psi(x):=0$ if $x$~is not a regular point.  

If $y$ is a regular value of~$\psi$, then 
the fiber $F_y:=\psi^{-1}(y)$ is a Riemannian submanifold 
of~$M$ of dimension $m-n$.   
Sard's lemma states that almost all $y\in N$ are regular values.

We can now state the {\em coarea formula}.

\begin{theorem}[Coarea formula]\label{pro:coarea}
Suppose that $M,N$ are Riemannian manifolds of 
dimensions~$m$, $n$, 
respectively, and let $\psi\colon M\to N$ be a surjective 
smooth map. Put $F_y=\psi^{-1}(y)$.
Then we have for any function $\chi\colon M\to\R$ that is integrable
with respect to the volume measure of $M$ that
\begin{equation}\tag*{\qed}
  \int_M \chi \, dM= \int_{y\in N} \left( \int_{F_y}
  \frac{\chi}{\NJ\psi}\, dF_y\right) dN .
\end{equation}
\end{theorem}

It should be clear that this result contains 
the transformation formula as a special case. 
Moreover, if we apply the coarea formula to the projection 
$\pi_2\colon M\times N \to N,\, (x,y)\mapsto y$, 
we retrieve Fubini's equality since $\NJ\pi_2=1$. 

The coarea formula is useful to define the concepts of marginal 
and conditional distributions for densities defined on a product 
space $M\times N$ when the components are Riemannian manifolds.

Suppose that we are in the situation described in the statement
of Theorem~\ref{pro:coarea} and we
have a probability measure on $M$ with density $\rho_M$.
For a regular value $y\in N$ we set
\begin{equation}\label{eq:pushf}
   \rho_N(y) := \int_{F_y} \frac{\rho_M}{\NJ\psi}\, dF_y .
\end{equation}
The coarea formula implies that for all measurable sets 
$B\subseteq N$ we have
$$
  \int_{\psi^{-1}(B)} \rho_M\, dM = \int_B \rho_N\, dN.
$$
Hence $\rho_N$ is a probability density on $N$.
We call it the {\em pushforward}
of $\rho_M$ with respect to $\psi$. 
In the special case that 
$\psi\colon M\times N\to N,\,(x,y)\mapsto y$,  
is  the projection, we have $\NJ\psi=1$, and 
we retrieve the usual formula for the marginal density.  

Furthermore, for a regular value $y\in N$ and $x\in F_y$ we define
the {\em conditional density} on $F_y$ 
\begin{equation}\label{eq:condd}
\rho_{F_y}(x) := \frac{\rho_M(x)}{\rho_N(y)\NJ\psi(x)} .
\end{equation}
Clearly, this defines a probability density on $F_y$.
Again, in the special case that 
$\psi\colon M\times N\to N,\,(x,y)\mapsto y$,  
we retrieve the usual formula for the conditional density.  

The coarea formula implies that for all measurable
functions $\chi\colon M\to \R$, 
$$
   \int_M \chi\, \rho_M \, dM= \int_{y\in N} \left(\int_{F_y} \chi\,
         \rho_{F_y}\, dF_y \right)\rho_N(y)\, dN ,
$$
provided the left-hand integral exists.
Therefore, we can interpret $\rho_{F_y}$ as the
{\em density of the conditional distribution} 
of $x$ on the fiber $F_y$ and
briefly express the formula above in probabilistic terms as
\begin{equation}\label{eq:Eit}
   \Exp_{x\sim \rho_M} \chi(x) = \Exp_{y\sim \rho_N}
   \Exp_{x\sim \rho_{F_y}}\chi(x).
\end{equation}

\subsection{An auxiliary result from linear algebra}
\label{se:aux-linalg} 

Let $E$ and $F$ be finite dimensional Euclidean vector spaces such 
that $\dim E\ge \dim F$. If $\varphi\colon E\to F$ is a surjective linear 
map, we denote by $\tilde{\varphi}\colon (\ker\varphi)^\perp \to F$ 
the restriction of $\varphi$ to the orthogonal complement of the 
kernel of $\varphi$. Then $\tilde{\varphi}$ is surjective and we define 
its {\em normal determinant} by 
$$
 \ndet\varphi := |\det\tilde{\varphi}| .
$$
(If $\varphi$ is not surjective, we set $\ndet\varphi := 0$.)
We consider the graph
$\Gamma:=\{(x,\varphi(x)) \mid x \in E\}$ of $\varphi$. Then, 
$\Gamma$ is a linear subspace of $E\times F$ and the two 
projections 
$$
  p_1\colon\Gamma \to E,\ (x,\varphi(x)) \mapsto x,\quad 
  p_2\colon\Gamma \to F,\ (x,\varphi(x)) \mapsto \varphi(x) 
$$
are linear maps.
Note that $p_1$ is an isomorphism and $p_2$ is surjective 
as $\varphi$ is so. 

\begin{lemma}\label{le:detquot}
Under the above assumptions, we have 
$$
   \frac{\ndet p_1}{\ndet p_2} = (\ndet\varphi)^{-1} .
$$
\end{lemma}

\proof
Let $K:=\ker\varphi$ and $\tilde{E}$ be the orthogonal complement of 
$K$ in $E$. 
Let $\tilde{\Gamma}\subseteq \tilde{E}\times F$ denote the 
graph of $\tilde{\varphi}\colon\tilde{E}\to F$. 
Further, let $\tilde{p}_1\colon\tilde{\Gamma}\to\tilde{E}$ and 
$\tilde{p}_2\colon\tilde{\Gamma}\to F$ denote the projections. 
Since $\varphi$ is surjective, $\tilde{\varphi}$ is bijective.  
We have obvious isometries
$\Gamma\simeq\tilde{\Gamma}\times K$, 
$p_1\simeq \tilde{p}_1\times \Id_K$, 
and we can interpret $\tilde{\Gamma}$ as the orthogonal complement 
of $\ker p_2 =K \times\{0\}$ in $\Gamma$. 
By the definition of the normal determinant, we have 
$$
 \ndet\varphi = |\det\tilde{\varphi} | , \quad
 \ndet p_1 = | \det \tilde{p}_1 |, \quad 
 \ndet p_2 = | \det \tilde{p}_2 | . 
$$
It is therefore sufficient to prove that 
$$
  \Big| \frac{\det\tilde{p}_1}{\det\tilde{p}_2} \Big| 
  = |\det\tilde{\varphi}|^{-1} .
$$
The singular value decomposition tells us that, with respect to suitable 
orthonormal bases on $\tilde{E}$ and $F$, 
the representation matrix of $\tilde{\varphi}$ equals 
$\diag(\s_1,\ldots,\s_n)$, where 
$\s_1\ge\cdots\ge\s_n$ are the singular values of $\tilde{\varphi}$. 
Note that $|\det\tilde{\varphi}|= \s_1\cdots \s_n$. 
It is now straightforward to check that 
$$
 |\det\tilde{p}_1| = \prod_{i=1}^n \frac{1}{\sqrt{1+\s_i^2}}, \quad 
 |\det\tilde{p}_2| = \prod_{i=1}^n \frac{\s_i}{\sqrt{1+\s_i^2}} .
$$
Therefore,
$$
 \Big| \frac{\det \tilde{p}_1}{\det \tilde{p}_2} \Big| 
 = \prod_{i=1}^n \sigma_i^{-1} = |\det\tilde{\varphi}|^{-1},
$$
which finishes the proof.
\eproof

Suppose that $W$ is a finite dimensional complex vector space 
with a Hermitian inner product $\langle~,~\rangle_\C$. 
The real part of $\langle~,~\rangle_\C$ turns $W$ into a 
Euclidean vector space. Let 
$\psi\colon W\to W$ be a $\C$-linear map.
If we denote by $\det\psi_\R$ the determinant of $\psi$, considered 
as an $\R$-linear map, then it is a well known fact that 
$\det\psi_\R = |\det\psi|^2$. 
(Indication of proof: the singular value decomposition allows to 
reduce to the case $W=\C$ and $\psi(z)=\lambda z$. 
Since $\psi$ is the composition of 
a rotation and a homothety by the stretching factor~$|\lambda|$, 
it follows that $\det\psi_\R = |\lambda|^2$.)

\subsection{Normal Jacobians for the eigenpair problem} 

Recall from \S\ref{subsec:spaces}, we have the two projections 
$$
 \pi_1\colon \V\to X, (A,\lambda,[v]) \to A, \quad
 \pi_2\colon \V\to Y, (A,\lambda,[v]) \to (\lambda,[v])
$$
and, for $(A,\lambda,v)\in\V$, the linear operator 
$A_{\lambda,v}\colon T_v \to T_v$ given by 
$P_{v^\perp}(A-\lambda\,\Id)_{|T_v}$.   

\begin{proposition}\label{pro:NJquot}
Let  $p:=(A,\lambda,[v]) \in \V$. Then $\lambda$ is a simple 
eigenvalue of $A$ iff 
$A_{\lambda,v}$ is invertible. In this case, the derivative
$D\pi_1(p)\colon T_p \to T_AX$ is an isomorphism, the derivative 
$D\pi_2(p)\colon T_p \to T_{(\lambda,[v])}Y$ is surjective, and we have 
$$
 \frac{\NJ\pi_1(p)}{\NJ\pi_2(p)} = | \det A_{\lambda,v} |^{2} 
 = \det(A_{\lambda,v}A_{\lambda,v}^*).
$$
\end{proposition}

\proof
Let $p:=(A,\lambda,[v]) \in \V$. We suppose that $\|v\|=1$ and we 
identify the tangent space $T_{[v]} \P(\C^n)$ with $T_v$. 
By orthogonal  invariance, we may assume without loss of generality 
that $v=e_1=(1,0,\ldots,0)$. Then we write 
$$
 A = \begin{bmatrix} \lambda & c^T \\ 0 & B \end{bmatrix}, \quad 
  c \in\C^{n-1}, B \in \C^{(n-1)\times (n-1)} .
$$ 
The matrix $M :=\lambda I - B$ represents the linear map 
$A_{\lambda,v}$. Clearly, $M$ is invertible iff $\lambda$ is a simple 
eigenvalue of $A$, which shows the first assertion of the proposition. 
We assume now that $M$ is invertible. 

Let $(\dot{A},\dot{\lambda},\dot{v})\in\C^{n\times n}\times\C\times T_v$.
According to equation (14.19) in~\cite{Condition}, 
the tangent space of $\V$ at $p$ is characterized in the following way:
\begin{equation}\label{eq:char-TV}
 (\dot{A},\dot{\lambda},\dot{v}) \in T_p\V \Longleftrightarrow 
 \dot{A}v + (A-\lambda I)\dot{v} - \dot{\lambda} v = 0. 
\end{equation}
In order to express $\dot{\lambda}$, $\dot{v}$ in terms of 
$\dot{A}$, we denote by $\dot{a}_i \in\C^n$ the $i$th column of $A$,  
and we write $\dot{a}_1 =(\dot{a}, \dot{b})$ where $\dot{a}\in\C$ and 
$\dot{b}\in\C^{n-1}$. Also, since $\dot{v}\in T_v$ (and we are assuming 
$v=e_1$) we have $\dot{v} = (0,\dot{w})$ for some $\dot{w}\in\C^{n-1}$.
Using this notation, equation~\eqref{eq:char-TV} can be rewritten as 
\begin{eqnarray*}
  \dot{a} + c^T\dot{w} - \dot{\lambda} & = & 0 \\
\dot{b}  - M\dot{w} &=& 0 .
\end{eqnarray*}
This system of equations has the unique solution
\begin{equation}\label{eq:sol}
  \dot{w} = M^{-1} \dot{b},\quad 
  \dot{\lambda} = \dot{a} + c^T M^{-1} \dot{b} .
\end{equation}
So we can interpret $T_p\V$ as the graph $\Gamma$ of the linear 
map $\varphi\colon E\to F$ (with 
$E=\C^{n\times n}$ and $F=\C\times T_v$) given by 
$\dot{A}\mapsto (\dot{\lambda},\dot{v})$,  
and $D\pi_1(p), D\pi_2(p)$ are the corresponding projections 
$p_1\colon\Gamma\to E$ and $p_2\colon\Gamma\to F$, respectively. 
According to Lemma~\ref{le:detquot}, it therefore suffices to prove that 
\begin{equation}\label{eq:ndetvarphi}
 \ndet\varphi = |\det A_{\lambda,v}|^{-2} .
\end{equation}

Since $\dot{w}=0$, $\dot{\lambda}=0$ implies $\dot{b} =0$, $\dot{a}=0$, 
we see that the orthogonal complement of the kernel of $\varphi$ is 
given by the conditions $\dot{a}_2=0,\ldots,\dot{a}_n=0$. The restriction 
$\tilde{\varphi}$ of $\varphi$ to $(\ker\varphi)^\perp$, 
$(\dot{a},\dot{b}) \mapsto (\dot{\lambda},\dot{w})$, 
according to~\eqref{eq:sol}, has the following matrix  
$$
 \begin{bmatrix} 1 & c^T M^{-1} \\ 0 & M^{-1} \end{bmatrix}
$$
with respect to the standard bases. Therefore, 
$\det\tilde{\varphi} = \det M^{-1}$. The determinant of $\tilde{\varphi}$, 
seen as a $\R$-bilinear map, therefore equals $|\det M^{-1}|^2$, 
see the comment at the end of \S\ref{se:aux-linalg}.  
We conclude that 
$$
 \ndet\varphi = |\det M |^{-2} = |\det A_{\lambda,v}|^{-2} ,
$$
which proves the claimed equality \eqref{eq:ndetvarphi}. 
\eproof

\subsection{Orthogonal decompositions}

We shall distinguish points in $\P(\C^n)$ from their representatives 
in $\IS(\C^n)$ and, accordingly, consider the following lifting 
of the solution variety $\V$
$$
    \aV:=\{(A,\lambda,v)\in\C^{n\times n}\times \C\times\IS(\C^n)
             \mid Av=\lambda v\}.
$$
Abusing notation, we denote the projection
$\aV\to\C\times \IS(\C^n)$ by $\pi_2$ as well. 
The fiber of $\pi_2$ at $(\lambda,v) \in \C\times\IS(\C^n)$ equals 
$$
  V_{(\lambda,v)}=\{A\in\Cnn\mid Av=\lambda v\}.
$$
This is an affine linear subspace of $\Cnn$ with the corresponding linear space 
$V^\lin_{(\lambda,v)} =\{ A\in\Cnn \mid Av=0\} $. 
We denote by $C_{(\lambda,v)}$ the orthogonal 
complement of $V^\lin_{(\lambda,v)}$ in $\Cnn$. 
So we have the orthogonal decomposition 
\begin{equation}\label{eq:odecomp2}
    \Cnn=V^\lin_{(\lambda,v)}\oplus C_{(\lambda,v)}. 
\end{equation} 
Let $\bH$ denote the point in $V_{(\lambda,v)}$ that is 
closest to the origin (with respect to the Frobenius norm). 
Note that $\bH\in C_{(\lambda,v)}$ and 
$V_{(\lambda,v)} = V^\lin_{(\lambda,v)} + K_{(\lambda,v)}$. 



Recall, $\varphi^{\oA,\sigma}_{n\times n}$ 
denotes the density of the Gaussian $\mcN(\oA,\s^2\Id)$ on $\Cnn$,
where $\oA\in\Cnn$ and $\s>0$.
For fixed $(\lambda,v)\in \aS$, we decompose the mean~$\oA$ 
acccording to~\eqref{eq:odecomp2} as 
$$
   \oA = \oM_{(\lambda,v)} +  \oK_{(\lambda,v)} 
$$
where $\oM_{(\lambda,v)}\in V^\lin_{(\lambda,v)}$ and 
$\oK_{(\lambda,v)}\in C_{(\lambda,v)}$. 
If we denote by $\varphi_{V^\lin_{(\lambda,v)}}$ 
and $\varphi_{C_{(\lambda,v)}}$ the densities of the
Gaussian distributions in the spaces 
$V^\lin_{(\lambda,v)}$ and $C_{(\lambda,v)}$ with 
covariance matrices $\s^2\Id$ and means 
$\oM_{(\lambda,v)}$ and $\oK_{(\lambda,v)}$,
respectively, then the density $\varphi^{\oA,\sigma}_{n\times n}$ 
factors as
\begin{equation}\label{eq:dfactor}
    \varphi^{\oA,\sigma}_{n\times n}(M+K) 
   = \varphi_{V^\lin_{(\lambda,v)}}(M)
\cdot \varphi_{C_{(\lambda,v)}}(K).
\end{equation}

Consider now the projection $\pi_1\colon{\aV}\to \Cnn$. 
Its fiber $V_A$ at $A\in\Cnn\setminus\Sigma$
is a disjoint union of $n$ unit circles and therefore has 
volume $2\pi\Dn$.
We associate with $\varphi^{\oA,\sigma}_{n\times n}$ the function 
$\rho_{\aV}:\aV\to\R$ defined by
\begin{equation}\label{eq:dens_V1}
  \rho_{\aV}(A,\lambda,v):=\frac{1}{2\pi n}\,
  \varphi^{\oA,\sigma}_{n\times n}(A)\,
  \NJ\pi_1(A,\lambda,v) .
\end{equation}
The proof of the following result is done 
as in~\cite[Lemma~18.10]{Condition}.

\begin{lemma}\label{le:rho1*}
{\bf (a)\ } 
The function $\rho_{\aV}$ is a probability density on $\aV$.
\smallskip

{\bf (b)\ }  
The expectation of a function
$F\colon\aV\to \R$ with respect to~$\rho_{\aV}$ can be expressed as
$$
   \Exp_{(A,\lambda,v)\sim \rho_{\aV}} F(A,\lambda,v) 
 = \Exp_{A\sim \varphi^{\oQ,\sigma}_{n\times n}} F_{\mathsf{sav}}(A),
$$ 
where
$F_{\mathsf{sav}}(A) :=  \frac1{2\pi n}\int_{V_A} F\, dV_A$.
\smallskip

{\bf (c)\ }  
The pushforward of $\rho_{\aV}$ with respect to 
$\pi_1\colon \aV\to\Cnn$ equals $\varphi^{\oA,\sigma}_{n\times n}$.
\smallskip

{\bf (d)\ }  
For $A\not\in\Sigma$, the conditional density on 
the fiber $V_A$  is the
density of the uniform distribution on $V_A$. \eproof
\end{lemma}

\begin{remark}\label{rem:average}
In the particular case that $F\colon\aV\to\R$ is given by 
$F(A,\lambda,v)=\frac{\mu^2(A,\lambda,v)}{\|A\|^2_F}$ 
we have $F_{\mathsf{sav}}(A)=\frac{\mu_{\av}^2(A)}{\|A\|^2_F}$. 
\end{remark}

Recall from \eqref{eq:defAlv} that we assigned to $A\in\C^{n\times n}$ 
the linear operator $A_{\lambda,v}\colon T_v\to T_v$ obtained by 
restricting $A-\lambda\Id$ to $T_v$. 
We write 
\begin{equation}\label{eq:defc}
 c_{(\lambda,v)} := \int_{M\in V^\lin_{(\lambda,v)}} |\det M_{\lambda,v}|^2\, \varphi_{V^\lin_{(\lambda,v)}}(M)\, dM .
\end{equation}

\begin{lemma}\label{pro:rhoW-2}
The pushforward density 
$\rho_{\aS}$ of $\rho_{\aV}$ with respect to
$\pi_2\colon \aV\to \aS$ is given by
\begin{equation}\label{eq:1st-stat}
  \rho_{\aS}(\lambda,v) = \frac{c_{(\lambda,v)}}{2\pi n}\cdot  
  \varphi_{C_{(\lambda,v)}}(\bH),
\end{equation}
the conditional density $\tilde{\rho}_{V_{(\lambda,v)}}$ 
on the fiber $V_{(\lambda,v)}$ of~$\pi_2$ is given by
\begin{equation*}
   \tilde{\rho}_{V_{(\lambda,v)}}(A) 
  = c_{(\lambda,v)}^{-1}\cdot \det (A_{\lambda,v}A_{\lambda,v}^*) 
    \varphi_{V^\lin_{(\lambda,v)}}(A-\bH),
\end{equation*}
and we have
$$ 
 \frac{\rho_{\aV}}{\NJ \pi_2}(A,\lambda,v) 
 = \rho_{\aS}(\lambda,v)\cdot \tilde{\rho}_{V_{(\lambda,v)}}(A).
$$
\end{lemma}

\proof
It is easy to check that Proposition~\ref{pro:NJquot} remains valid 
when replacing $\P(\C^n)$ by $\IS(\C^n)$. So for $(A,\lambda,v)\in\aV$ we get 
\begin{equation}\label{eq:NJqneu}
 \frac{\NJ \pi_1}{\NJ \pi_2}(A,\lambda,v) = 
 |\det (A_{\lambda,v}) |^2 
\end{equation}
We write $A=M + \bH$ with $M \in V^\lin_{(\lambda,v)}$ .
Combining \eqref{eq:NJqneu} with~\eqref{eq:dens_V1}, we get
\begin{eqnarray}\label{eq:primera}
 \frac{\rho_{\aV}}{\NJ \pi_2}(A,\lambda,v)
  &=&\frac1{2\pi n}\, \varphi^{\oQ,\sigma}_{n\times n}(A)\cdot 
   |\det(A_{\lambda,v})|^2\\ 
  &=& \frac1{2\pi n}\, \varphi_{V^\lin_{(\lambda,v)}}(M)\cdot 
   \varphi_{C_{(\lambda,v)}}(\bH)\cdot 
   |\det(A_{\lambda,v})|^2 . 
\end{eqnarray}
For fixed $(\lambda,v)$ we integrate both sides of this equality 
over $M\in V^\lin_{(\lambda,v)}$. Equality~\eqref{eq:pushf} tells us that 
on the left-hand side we obtain $\rho_{\aS}(\lambda,v)$. On the 
right-hand side we obtain 
$\frac1{2\pi n}\,c_{(\lambda,v)}\cdot \varphi_{C_{(\lambda,v)}}(\bH)$ 
by the definition~\eqref{eq:defc} of $c_{(\lambda,v)}$.
This proves the first equality~\eqref{eq:1st-stat}
in the statement. 

For the second, we use the definition~\eqref{eq:condd} for the 
conditional density and 
\eqref{eq:dens_V1}, 
\eqref{eq:1st-stat},
\eqref{eq:dfactor}, and 
\eqref{eq:NJqneu} to get 
\begin{eqnarray*}
 \tilde{\rho}_{V_{(\lambda,v)}}(A) &=& \frac{\rho_{\aV}(A,\lambda,v)}
        {\rho_{\aS}(\lambda,v)\,\NJ\pi_1(A,\lambda,v)}\notag\\
  &=& \frac{\frac{1}{2\pi n}\, \varphi^{\oA,\sigma}_{n\times n}(A)\, \NJ\pi_1(A,\lambda,v)}{
                       {\frac{1}{2\pi n}\, c_{(\lambda,v)} \,\varphi_{C_{(\lambda,v)}}(\bH)\, \NJ\pi_2(A,\lambda,v)} }\\
  &=& c_{(\lambda,v)}^{-1}\cdot \varphi_{V^\lin_{(\lambda,v)}}(M)\cdot |\det(A_{\lambda,v})|^2 . 
\end{eqnarray*}
This proves the second inequality. 
The third equality is a trivial consequence 
of~\eqref{eq:primera}
and the first two assertions of the lemma. 
\eproof

We can now give the proof of the main result in this section. 
\medskip

\proofof{Theorem~\ref{th:mu2-bound}}
Because of Lemma~\ref{le:rho1*}(b) (and Remark~\ref{rem:average}) 
we have
\begin{equation*}\label{eq:HV}
  \Exp_{Q\sim \mcN(\oQ,\s^2\Id)}
  \Big(\frac{\mu_{\av}^2(Q)}{\|Q\|_F^2} \Big) 
  \ =\
  \Exp_{(Q,\lambda,v)\sim\rho_{\aV}} \Big(\frac{\mu^2(Q,\lambda,v)}
  {\|Q\|_F^2} \Big).
\end{equation*}
By the definition of the condition number~\ref{eq:defmu} we have 
\begin{equation*}
   \frac{\mu(Q,\lambda,v)}{\|Q\|_F} 
    \,=\,\|Q_{\lambda,v}^{-1}\|.
\end{equation*}
Hence 
\begin{eqnarray}\label{eq:EVEM}
 \Exp_{(Q,\lambda,v)\sim\rho_{\aV}} 
 \Big(\frac{\mu^2(Q,\lambda,v)}{\|Q\|_F^2} \Big) &=&
 \Exp_{(Q,\lambda,v)\sim\rho_{\aV}}
 \big(\|Q_{\lambda,v}^{-1}\|^2\big)\notag \\
&=& \Exp_{(\lambda,v)\sim\rho_{\aS}}
  \Big( \Exp_{Q\sim \tilde{\rho}_{V_{(\lambda,v)}}}
  \big( \|Q_{\lambda,v}^{-1}\|^2\big)  \Big)
\end{eqnarray}
the last by equation~\eqref{eq:Eit}. 

Because of unitary invariance, we may assume $v=e_1=(1,0,\ldots,0)$. In this case, 
we have 
$$
  V_{(\lambda,e_1)} = \Big\{
  \begin{bmatrix}  \lambda & a \\ 0 & B\end{bmatrix} \ \Big| \ 
    a\in\C^{n-1}, B\in\C^{(n-1)\times(n-1)} \Big\}, \quad 
  \bH = \begin{bmatrix}  \lambda & 0 \\ 0 & 0\end{bmatrix}, 
$$
and $V^\lin_{(\lambda,e_1)}$ equals the space of matrices of the form 
$M=\begin{bmatrix}  0 & a \\ 0 & B\end{bmatrix}$. 
The Gaussian distribution 
$\varphi_{V^\lin_{(\lambda,v)}}$ of $M$ 
induces an isotropic Gaussian distribution
$\rho_\sigma^{(n-1)}$ of $a\in\C^{(n-1)}$ with center $\oa$ 
and an isotropic Gaussian distribution
$\rho_\sigma^{(n-1)^2}$
of $B\in \C^{(n-1)\times(n-1)}$ with center $\oB$.

Decomposing $Q\in V_{(\lambda,e_1)}$ as $Q=M + \bH$, 
Lemma~\ref{pro:rhoW-2} tells us that 
the conditional density $\tilde{\rho}_{V_{(\lambda,e_1)}}$ 
on the fiber $V_{(\lambda,e_1)}$ has the form 
\begin{eqnarray}\label{eq:form}
\tilde{\rho}_{V_{(\lambda,e_1)}}(Q) &=& 
c_{(\lambda,e_1)}^{-1}\cdot |\det(Q_{\lambda,e_1})|^2 
     \rho_{V^\lin_{(\lambda,e_1)}}(M)\notag\\
&=&c_{(\lambda,e_1)}^{-1}\cdot |\det(B-\lambda\Id)|^2\, 
     \rho^{(n-1)}_\sigma(a)\, \rho^{(n-1)^2}_\sigma(B).
\end{eqnarray}
For the second equality we used that $Q_{\lambda,v}$ is represented by the matrix 
$B- \lambda\Id$, see~\eqref{eq:defAlv}. 
By the definition~\eqref{eq:defc} of $c_{(\lambda,e_1)}$, 
\begin{eqnarray}\label{eq:c_lambda}
c_{(\lambda,e_1)} &=& 
   \Exp_{M\sim\rho_{V^\lin_{(\lambda,e_1)}}} 
   |\det(M_{\lambda,v})|^2\notag\\ 
&=&\int_{a\in\C^{(n-1)}\atop B\in\C^{(n-1)\times(n-1)}}
        |\det(B-\lambda\Id)|^2\,
        \rho^{(n-1)}_\sigma(a)\, \rho^{(n-1)^2}_\sigma(B)\, da
       \,dB\notag\\
&=&\Exp_{\tilde{B}\sim \mcN(\oB-\lambda\Id,\sigma^2\Id)} |\det \tilde{B}|^2.
\end{eqnarray}
It follows from~\eqref{eq:form} 
that,
\begin{eqnarray*}
  \Exp_{Q\sim \tilde{\rho}_{V_{(\lambda,e_1)}}}
  \big( \|Q_{\lambda,e_1}^{-1}\|^2\big) 
& =&\int_{a\in\C^{(n-1)}\atop B\in\C^{(n-1)\times(n-1)}}
     \|(B-\lambda\Id)^{-1}\|^2\;
    c_{(\lambda,e_1)}^{-1}\, 
  |\det(B-\lambda\Id)|^2\\
  & &\qquad \rho^{(n-1)}_\sigma(a)\, 
  \rho^{(n-1)^2}_\sigma(B)\,da\,dB\\
&=& \Exp_{B\sim \mcN(\oB,\sigma^2\Id)} 
     \|(B-\lambda\Id)^{-1}\|^2\, c_{(\lambda,e_1)}^{-1}\, 
    |\det(B-\lambda\Id)|^2\\
&=& \Exp_{\tilde{B}\sim \mcN(\oB-\lambda\Id,\sigma^2\Id)} 
\|\tilde{B}^{-1}\|^2\, c_{(\lambda,e_1)}^{-1}\,|\det \tilde{B}|^2.
\end{eqnarray*}
The form of this expectation (along with that of 
$c_{(\lambda,e_1)}$ given in~\eqref{eq:c_lambda}) is 
exactly the one in the hypothesis 
of~\cite[Proposition~4.22]{Condition}, a result then 
ensuring that
$$
\Exp_{Q\sim \tilde{\rho}_{V_{(\lambda,e_1)}}}
  \big( \|Q_{\lambda,e_1}^{-1}\|^2\big) \leq \frac{en}{2\sigma^2}.
$$
The fact that this last bound does not depend on the pair 
$(\lambda,e_1)$ 
allows to quickly derive our result from~\eqref{eq:EVEM}.
\eproof

\section{Proof of Theorems~\ref{thm:main}
and~\ref{thm:main2}}\label{sec:main_proof}

\subsection{A useful change of variables}\label{subsec:change}

Since a linear combination (with fixed coefficients)
of two Gaussian matrices is Gaussian as well, it is convenient
to parameterize the interval $[M,A]$ by a parameter $t\in[0,1]$
representing a ratio of Euclidean distances (instead of a ratio 
of angles as $\tau$ does).
Thus we write, abusing notation, $Q_t=tA+(1-t)M$.
For fixed $t$, as noted before, $Q_t$ follows a Gaussian law. 
For this new parametrization we have the following result 
(see~\cite[Proposition~5.2]{BuCu11} for a proof).

\begin{proposition}\label{prop:a-ALH}
Let $A,M\in\Cnn$ be $\R$-linearly independent and 
$\alpha:=\dS(A,M)$. The function
\begin{eqnarray*}
    t:[0,1]&\to&[0,1]\\
        \tau&\mapsto&t(\tau):=\frac{\|M\|_F}
     {\|M\|_F + \|A\|_F (\sin\a\cot(\tau\a)-\cos\a)}
\end{eqnarray*}
is a bijection satisfying, for every $\tau\in[0,1]$, that 
$$
   Q_\tau=t(\tau)A+(1-t(\tau))M.
$$
Furthermore, for all $0\leq a\leq b\leq 1$, 
\begin{equation}\tag*{\qed}
 \dS(A,M)\int_{a}^b \mu_{\av}^2 (Q_\tau)d\tau
 \leq \|A\|_F\,\|M\|_F\,
  \int_{t(a)}^{t(b)} \frac{\mu_{\av}^2(Q_t)}{\|Q_t\|_F^2}\, dt.
\end{equation}
\end{proposition}

\subsection{Proof of Theorem~\ref{thm:main}}

We want to bound $\avcost(n)=\Oh(n^3)\aviter(n)$. 
To do so it is enough to bound $\aviter(n)$. Recall, we have 
\begin{eqnarray}\label{eq:step1}
  \aviter(n)&=& \E_{A\sim\mcN(0,\Id)}\frac1n
        \sum_{j=1}^n K(A,M,m_{jj},e_j)\notag\\
  &\le& 1077\E_{A\sim\mcN(0,\Id)}
            \dS(M,A)\frac1n\sum_{j=1}^n \int_0^1
          \mu^2(Q_\tau,\lambda_\tau^{(j)},v_\tau^{(j)})d\tau\notag\\     
  &=& 1077\,\int_0^1\E_{A\sim\mcN(0,\Id)}
            \dS(M,A)\mu_{\av}^2(Q_\tau)d\tau\\     
  &=& 1077\,\int_0^1\E_{A\sim\mcN_T(0,\Id)}
            \dS(M,A) 
            \mu_{\av}^2(Q_\tau)d\tau\notag
\end{eqnarray}
the second line by Proposition~\ref{thm:main_path_following}, 
the third by the definition of $\mu_{\av}$, and  
the fourth because, for each $\tau\in[0,1]$, the function 
$A\mapsto\dS(M,A)\mu_{\av}^2(Q_\tau)$ is scale invariant 
(and we therefore use~\eqref{eq:truncating}). 

We are left with the task of bounding the last expression 
in~\eqref{eq:step1}. The general idea is to use the change of 
variables described in~\S\ref{subsec:change}, then use 
Theorem~\ref{th:mu2-bound} to bound the inner expectation, 
and finally integrate over $t\in[0,1]$. This direct approach is
however infeasible, because the resulting integral in $t$ 
would be improper. To circumvent this difficulty, the idea 
(going back to~\cite{BuCu11}) is simple. For small values of~$\tau$ 
the matrix $Q_\tau$ is close to $M$, and therefore,
the value of $\mu_{\av}^2(Q_\tau)$ can be bounded by
a small multiple of $\mum^2(M)$. For the remaining values of 
$\tau$, the corresponding $t=t(\tau)$ is bounded away 
from 0, and therefore, so is the variance 
in the distribution $\mcN(\oQ_t,t^2\Id)$
for $Q_t$ (here $\oQ_t:=(1-t)M$). 
This allows one to control the denominator on 
the right-hand side of Theorem~\ref{th:mu2-bound} 
when using this result. Here are the details.

We set $\alpha:=\dS(A,M)$, $\e:=0.12$, 
$C_\e:=\frac{\e}{12.5}=0.0096$, and  
$\xi:=\frac{2C_\e(1-\e)}{3\sqrt{6}(1+\e)^4}\approx 0.001461$, 
as in the proof of Theorem~\ref{thm:main_path_following}, 
and define 
$$
  T := \sqrt{2}\, n,\quad 
   \delta_0 := \frac{\xi}{\mum^2(M)} \qquad\mbox{and}\qquad
   t_T:=\frac{1}{1+T+1.0000015\,\frac{T}{\delta_0}}.
$$
Let $(\lambda^{(1)},v^{(1)}),\ldots,(\lambda^{(n)},v^{(n)})$ be 
the eigenpairs of $M$
and denote by
$(Q_\tau,\lambda_\tau^{(j)},v_\tau^{(j)})_{\tau\in[0,1]}$ the
lifting of $[M,A]$ in $\V$
corresponding to the initial triple $(M,\lambda^{(j)},v^{(j)})$.

Corollary~\ref{cor:b1} for $i=0$ implies the following:
for all $j$ and all
$\tau\leq \frac{\xi}{\a\,\mu^2(M,\lambda^{(j)},v^{(j)})}$
we have
$$
   \mu(Q_\tau,\lambda_\tau^{(j)},v_\tau^{(j)}) 
   \leq (1+\e)\mu(M,\lambda^{(j)},v^{(j)})
    \leq (1+\e)\mum(M) .
$$
In particular,
this inequality holds for all $j$ and all
$\tau\leq \frac{\delta_0}{\a}$, and hence for all such $\tau$, 
we have
\begin{equation}\label{eq:tau0}
   \mu_{\av}(Q_\tau) \leq (1+\e) \mum(M).
\end{equation}
Splitting the integral in the last last expression 
in~\eqref{eq:step1} at 
$\tau_0:=\min\big\{1,\frac{\delta_0}{\dS(A,M)}\big\}$, we obtain
\begin{eqnarray}\label{eq:step3}
   \int_0^1\E_{A\sim\mcN_T(0,\Id)}
            \alpha 
            \mu_{\av}^2(Q_\tau)d\tau 
   &=& \Exp_{A\sim N_T(0,\Id)}\Big(\alpha
        \int_0^{\tau_0} \mu_{\av}^2 (Q_\tau)\,d\tau\Big)\notag\\ 
    & &  + \Exp_{A\sim N_T(0,\Id)}\Big(\alpha
       \int_{\tau_0}^1 \mu_{\av}^2 (Q_\tau)\,d\tau\Big).
\end{eqnarray}
Using~\eqref{eq:tau0} we bound the first term on the 
right-hand side as follows: 
$$
  \Exp_{A\sim N_T(0,\Id)}\Big(\alpha
       \int_0^{\tau_0} \mu_{\av}^2 (Q_\tau)\,d\tau\Big)
    \leq \delta_0\, (1+\e)^2\, \mum^2(M) = (1+\e)^2\xi 
  \le 0.002.
$$
For bounding the second term, we assume without loss of generality 
that $\tau_0\le 1$. We then use 
Proposition~\ref{prop:a-ALH} to obtain that for a fixed $A$ (recall $\|M\|_F=1$) 
\begin{equation}\label{eq:step2.5}
   \alpha
       \int_{\tau_0}^1 \mu_{\av}^2 (Q_\tau)\,d\tau
  \leq \int_{t_0}^1 \|A\|_F \frac{\mu_{\av}^2 (Q_t)}{\|Q_t\|_F^2}\,dt
\end{equation}
with 
$$
   t_0=\frac{1} {1+\|A\|(\sin\a\cot\delta_0-\cos\a)}.
$$
Now note that $\|A\|_F\leq T$, since we draw $A$ from $N_T(0,\Id)$.
This allows us to bound~$t_0$ from below by a quantity 
independent of~$A$. Indeed, we first note that 
$$
    0 \le \sin\a \cot\delta_0 -\cos\a \le \frac{1}{\sin\delta_0} - \cos\a
    \le \frac{1}{\sin\delta_0} +1 ,
$$
and moreover,
$\sin\delta_0 \ge 0.9999985\cdot\delta_0$, 
since $\delta_0\le 2\xi \le 0.002922$ (Lemma~\ref{le:lb_mu}). 
We can now use that $\|A\|_F\leq T$ and bound $t_0$ as
$$
  t_0 \geq \frac{1}{1+T+\frac{T}{\sin \delta_0}}
  \geq \frac{1}{1+T+1.0000015\,\frac{T}{\delta_0}}=t_T.
$$
We next use this bound, together with~\eqref{eq:step2.5},
and bound  the second term in~\eqref{eq:step3}: 
\begin{align*}
\Exp_{A\sim N_T(0,\Id)}\Big(\alpha
      \int_{\tau_0}^1 &\mu_{\av}^2 (Q_\tau)\,d\tau\Big)
\,\leq\,
\Exp_{A\sim N_T(0,\Id)} \Big( T \int_{t_T}^1
         \frac{\mu_{\av}^2(Q_t)}{\|Q_t\|_F^2}\, dt \Big)\\
   =\;&    T \int_{t_T}^1
  \E_{A\sim N_T(0,\Id)}
        \bigg(\frac{\mu_{\av}^2(Q_t)}{\|Q_t\|_F^2}\bigg)\, dt
    \,\leq\,  \frac{T}{P_{T,1}} \int_{t_T}^1
    \Exp_{A\sim \mcN(0,\Id)}
    \bigg(\frac{\mu_{\av}^2(Q_t)}{\|Q_t\|_F^2}\bigg)\, dt.
\end{align*}
Observing that for fixed~$t$ and when $A$ is distributed
following $\mcN(0,\Id)$, the variable
$Q_t = (1-t)M + tA$ follows the Gaussian
$\mcN(\oQ_t, t^2\Id)$, we deduce 
(recall $T=\sqrt{2} n$ and $P_{T,1} \ge 1/2$ by Lemma~\ref{lem:X})  
\begin{equation*}
\int_0^1\E_{A\sim\mcN_T(0,\Id)}
            \a \mu_{\av}^2(Q_\tau)d\tau 
\leq 0.002
  + 2\sqrt{2}\,n\int_{t_T}^1
    \Exp_{Q_t\sim \mcN(\oQ_t,t^2\Id)}
    \bigg(\frac{\mu_{\av}^2(Q_t)}{\|Q_t\|_F^2}\bigg)\, dt.
\end{equation*}
To bound the integral in the right-hand side 
we apply Theorem~\ref{th:mu2-bound} and obtain
\begin{equation*}
\begin{split}
    \int_{t_T}^1
    \Exp_{Q_t\sim \mcN(\oQ_t,t^2\Id)}
    \Big(\frac{\mu_{\av}^2(Q_t)}{\|Q_t\|_F^2}\Big)\, dt
  & \;\leq\ \int_{t_T}^1\frac{en}{2t^2}\,dt
  =  \frac{en}{2}\bigg(\frac1{t_T}-1\bigg) 
  \;=\  \frac{enT}{2}\Big(1+\frac{1.0000015}{\delta_0} \Big) \\
  &\;=\  \frac{en^2}{\sqrt{2}}\Big(1+\frac{1.0000015\,\mum^2(M)}{\xi} \Big) 
  \;=\  \Oh(n^4)
\end{split}
\end{equation*}
the last by Lemma~\ref{lem:initial}. 
We conclude that 
$$
\Exp_{A\sim N_T(0,\Id)}\Big(\alpha
      \int_{\tau_0}^1 \mu_{\av}^2 (Q_\tau)\,d\tau\Big)
\,\leq\, 0.002 +\Oh(n^5)
$$
and hence, that $\aviter(n)=\Oh(n^5)$, and that 
$\avcost(n)=\Oh(n^8)$.
\medskip

We finally prove the smoothed analysis bounds. Reasoning as 
in~\eqref{eq:step1} we see that the smoothed number of iterations 
$\siter(n,\sigma)$ of {\sf Single\_Eigenpair} satisfies
\begin{equation}\label{eq:step4}
  \siter(n,\sigma) \ \le\  1077\,\sup_{\oA\in\IS(\Cnn)}\,\int_0^1
            \E_{A\sim\mcN_T(\oA,\sigma^2\Id)}
            \alpha\,
            \mu_{\av}^2(Q_\tau)d\tau.
\end{equation}
We deal with the integral as above to obtain
$$
\int_0^1\E_{A\sim\mcN_T(\oA,\sigma^2\Id)}
            \a \mu_{\av}^2(Q_\tau)d\tau 
\leq 0.002
  + 2(\sqrt{2}n+1)\int_{t_T}^1
    \Exp_{Q_t\sim \mcN(\oQ_t,\sigma^2t^2\Id)}
    \bigg(\frac{\mu_{\av}^2(Q_t)}{\|Q_t\|_F^2}\bigg)\, dt
$$
where we used $\|A\|_F\leq T+\|\oA\|_F = \sqrt{2}n+1$ and 
we now have $\oQ:=(1-t)M+t\oA$. The rest of the reasoning 
follows by noting that the term $\frac{1}{\sigma^2}$ can be factored 
out the integral and that the resulting bound for the 
integral in~\eqref{eq:step4}, $\Oh(\frac{n^5}{\sigma^2})$, is 
independent of $\oA$.
\eproof

\subsection{Proof of Theorem~\ref{thm:main2}}

It is immediate from the fact that, for any $A\in\Cnn$, the 
number of iterations performed by {\sf All\_Eigenpairs}---to compute the $n$ eigenpairs of $A$--- is $n$ times 
the (expected) number of iterations performed by 
{\sf Single\_Eigenpair}---to compute one such eigenpair.
\eproof

{\small

}

\end{document}